\documentclass [11pt]{article}
\usepackage{latexsym}
\usepackage{amssymb}
\usepackage{citesort}

\newcommand{\beq}{\begin{equation}}
\newcommand{\eeq}{\end{equation}}
\newcommand{\bea}{\begin{eqnarray}}
\newcommand{\eea}{\end{eqnarray}}
\newcommand{\barr}[1]{\begin{array}}
\newcommand{\earr}{\end{array}}
\newtheorem{theorem}{Theorem}[section]
\newtheorem{definition}{Definition}[section]
\newtheorem{proposition}{Proposition}[section]
\newtheorem{lemma}{Lemma}[section]
\newcommand{\bl}{\begin{lemma}}
\newcommand{\el}{\end{lemma}}
\newcommand{\bdf}{\begin{definition}}
\newcommand{\edf}{\end{definition}}
\newcommand{\bth}{\begin{theorem}}
\newcommand{\enth}{\end{theorem}}

\newcounter{abc}
\newenvironment{zoznamrom}{\setcounter{abc}{0}\begin{list}{\roman{abc})}
                       {\usecounter{abc}}}{\end{list}}

\newcounter{zozal}

\newcounter{znum}
\newenvironment{zoznamnum}{\setcounter{znum}{0}\begin{list}{\arabic{znum}.}
                       {\usecounter{znum}}}{\end{list}}

\def\kp{k_{+}} 
\def\km{k_{-}} 
\def\kz{k_{0}}
\def\mR{\mathbb R}
\def\mC{\mathbb C}
\def\mZ{\mathbb Z}
\def\pq{(p,q)}
\def\zl{\left(}
\def\zr{\right)}
\def\c+{\rlap{\ \raisebox{.2ex}{\scriptsize+}}\supset}

\makeatletter
\def\invddots{\mathinner{\mkern1mu\raise\p@\vbox{\kern7\p@\hbox{.}}\mkern2mu
\raise4\p@\hbox{.}\mkern2mu\raise7\p@\hbox{.}\mkern1mu}}
\makeatother

\title{Maximal Abelian Subalgebras of e(p,q) algebras}
\author {Z. Thomova and P. Winternitz}
\date{}

\begin{document}
\bibliographystyle{unsrt}
\maketitle
\centerline{\bf CRM -2516}

%%%%%%%%%%%%%%%%%%%%%%%%%%%%%%%%%%%%%%%%%%%%%%%%%%%%%%%%%%%%%%%%%%%
%%%%%%%%%%%%%%%%%%%%%%%%%%%%%%%%%%%%%%%%%%%%%%%%%%%%%%%%%%%%%%%%%%%
%%%%%%%%%%%%%%%%%%%%%                         %%%%%%%%%%%%%%%%%%%%%
%%%%%%%%%%%%%%%%%%%%%%%%%%%%%%%%%%%%%%%%%%%%%%%%%%%%%%%%%%%%%%%%%%%
%%%%%%%%%%%%%%%%%%%%%%%%%%%%%%%%%%%%%%%%%%%%%%%%%%%%%%%%%%%%%%%%%%%
\bigskip
\bigskip

\begin{abstract}
Maximal abelian subalgebras of one of the classical real inhomogeneous Lie 
algebras are constructed, namely those of the pseudoeuclidean Lie algebra $e\pq$. 
Use is made of the semidirect sum structure of $e\pq$ with the translations 
$T(p+q)$ as an abelian ideal. We first construct splitting MASAs that are 
themselves direct sums of abelian subalgebras of $o\pq$ and of subalgebras 
of $T(p+q)$. The splitting subalgebras are used to construct the complementary 
nonsplitting ones. We present general decomposition theorems and construct 
indecomposable MASAs for all algebras $e\pq$, $p \geq q \geq 0$. The case 
of $q=0$ and $1$ were treated earlier in a physical context. The case $q=2$ 
is analyzed here in detail as an illustration of the general results.

\bigskip

Les sous-alg\`ebres maximales ab\'eliennes (SAMAs) d'une alg\`ebre r\'eelle classique 
non-homog\`ene sont construites, en particulier, celles d'alg\`ebre de Lie 
pseudo-euclidienne $e\pq$. 
On utilise la structure de la somme semi-directe de $e\pq$ avec les translations 
$T(p+q)$ qui repr\'esente un id\'eal ab\'elien. Nous avons construit, en premier, 
les SAMAs "splitting", qui sont des sommes directes des sous-alg\`ebres ab\'eliennes 
de $o\pq$ et de sous-alg\`ebres de $T(p+q)$. Les sous-alg\`ebres ``splitting'' sont 
utilis\'ees pour  construire les sous-alg\`ebres complementaire -"nonsplitting". 
Nous pr\'esentons les th\'eor\`emes g\'en\'eraux de d\'ecomposition et nous construisons 
les SAMAs ind\'ecomposables pour toutes les alg\`ebres $e\pq$, $p \geq q \geq 0$. 
Les cas de $q=0$ et $1$ sont d\'ej\`a  trait\'es dans un context physique. 
Le cas $q=2$ est analys\'e ici en d\'etail comme une illustration des r\'esultats 
g\'en\'eraux.

\end{abstract}

\newpage

%%%%%%%%%%%%%%%%%%%%%%%%%%%%%%%%%%%%%%%%%%%%%%%%%%%%%%%%%%%%%%%%%%%
%%%%%%%%%%%%%%%%%%%%%%%%%%%%%%%%%%%%%%%%%%%%%%%%%%%%%%%%%%%%%%%%%%%
%%%%%%%%%%%%%%%%%%%%%                         %%%%%%%%%%%%%%%%%%%%%
%%%%%%%%%%%%%%%%%%%%%%%%%%%%%%%%%%%%%%%%%%%%%%%%%%%%%%%%%%%%%%%%%%%
%%%%%%%%%%%%%%%%%%%%%%%%%%%%%%%%%%%%%%%%%%%%%%%%%%%%%%%%%%%%%%%%%%%
%%%%%%%%%%%%%%%%%%%%%%%%%%%%%%%%%%%%%%%%%%%%%%%%%%%%%%%%%%%%%%%%%%%
%%%%%%%%%%%%%%%%%%%%%%%%%%%%%%%%%%%%%%%%%%%%%%%%%%%%%%%%%%%%%%%%%%%
%%%%%%%%%%%%%%%%%%%%%                         %%%%%%%%%%%%%%%%%%%%%
%%%%%%%%%%%%%%%%%%%%%%%%%%%%%%%%%%%%%%%%%%%%%%%%%%%%%%%%%%%%%%%%%%%
%%%%%%%%%%%%%%%%%%%%%%%%%%%%%%%%%%%%%%%%%%%%%%%%%%%%%%%%%%%%%%%%%%%

\section {Introduction}

The purpose of this article is to present a classification of the maximal 
abelian subalgebras (MASAs) of the pseudoeuclidean Lie algebra $e\pq$. Since this 
Lie algebra can be represented by a specific type of real matrices of dimension 
$(p+q+1) \linebreak \times (p+q+1)$, the subject of this article is placed squarely within 
a classical problem of linear algebra, the construction of sets of commuting 
matrices.

Most of the early papers in this direction \cite{fro,sch,krav} as well as more 
recent ones \cite{ger,cou,gus,tau,laff}, were devoted to commuting matrices within 
the set of all matrices of a given dimension. In other words, they studied abelian
subalgebras of the Lie algebras $gl(n,\mC)$ and $gl(n,\mR)$. For a historical 
review with many references see the book by Suprunenko and Tyshkevich \cite{st}.

Maltsev constructed all maximal abelian subalgebras of maximal dimension for 
all complex finite-dimensional simple Lie algebras \cite{mal}. An important 
subclass of MASAs are Cartan subalgebras, {\it i.e.} self-normalizing MASAs 
\cite{jac}. The simple complex Lie algebras, as well as the compact ones, 
have just one conjugacy class of Cartan subalgebras. The real noncompact forms 
of the simple Lie algebras can have several conjugacy 
classes of them. They have been classified by Kostant \cite{kos} and 
Sugiura \cite{sug}.

This article is part of a series, the aim of which is to construct all MASAs 
of the classical Lie algebras. Earlier articles were devoted to the classical 
simple Lie algebras, such as $sp(2n,\mR)$ and $sp(2n,\mC)$ \cite{sp_C}, 
$su\pq$ \cite{su_pq}, $o(n,\mC)$ \cite{verc} and $o\pq$ \cite{verop}. General 
results for MASAs of classical simple Lie algebras are presented in \cite{decomp}.
More recently MASAs of some inhomogeneous classical Lie algebras were studied, 
namely those of $e(n,\mC)$ \cite{kw}, $e(p,0)$ and $e(p,1)$ \cite{conf}. Here 
we consider $e\pq$ for all $p \geq q \geq 0$. The two special cases, $q=0$ and 
$q=1$, treated earlier, are of particular importance in physics and are also 
much simpler than the general case.

The motivation for a study of MASAs was discussed in previous articles \cite{sp_C,su_pq,verc,verop,decomp,kw,conf}. As a mathematical problem the 
classification of MASAs 
is an 
extension of the classification of individual elements of Lie algebras into conjugacy
classes \cite{bc,mal2,dj}. A classification of MASAs of classical Lie algebras is 
an important 
ingredient in the classification of all subalgebras of these algebras.

In applications in the theory of partial differential equations, MASAs provide 
coordinate systems in which invariant equations allow the separation of variables. 
More specifically, they provide "ignorable variables" not figuring in the 
corresponding metric tensors, when considering Laplace-Beltrami or Hamilton-Jacobi 
equations. In quantum physics they provide complete sets of commuting operators. 
In classical physics they provide integrals of motion in involution.

The classification problem is formulated in Section 2, where we also present 
some necessary definitions and explain the classification strategy. Section 3 
contains a brief summary of the known results on MASAs of $o\pq$ \cite{verop}. 
They are needed in the rest of this article and we reproduce them in a 
condensed 
form to make the article self-contained. Section 4 is devoted to splitting 
subalgebras of $e\pq$, {\it i.e.} subalgebras that are direct sums of subalgebras 
of the algebra 
$o\pq$ and those of the translation algebra $T(p+q)$. The complementary case of nonsplitting 
MASAs of $e\pq$ is the subject of Section 5. The results on MASAs of $e\pq$ obtained 
in Sections 4 and 5 are reformulated in terms of a decomposition of the underlying 
linear space $S\pq$ in Section 6. Indecomposable MASAs of $e\pq$ are described 
in the same section. Section 7 is devoted to a special case in which all 
results 
are entirely explicit, namely MASAs of $e(p,2)$.

%%%%%%%%%%%%%%%%%%%%%%%%%%%%%%%%%%%%%%%%%%%%%%%%%%%%%%%%%%%%%%%%%%%
%%%%%%%%%%%%%%%%%%%%%%%%%%%%%%%%%%%%%%%%%%%%%%%%%%%%%%%%%%%%%%%%%%%
%%%%%%%%%%%%%%%%%%%%%                         %%%%%%%%%%%%%%%%%%%%%
%%%%%%%%%%%%%%%%%%%%%%%%%%%%%%%%%%%%%%%%%%%%%%%%%%%%%%%%%%%%%%%%%%%
%%%%%%%%%%%%%%%%%%%%%%%%%%%%%%%%%%%%%%%%%%%%%%%%%%%%%%%%%%%%%%%%%%%

%%%%%%%%%%%%%%%%%%%%%%%%%%%%%%%%%%%%%%%%%%%%%%%%%%%%%%%%%%%%%%%%%%%
%%%%%%%%%%%%%%%%%%%%%%%%%%%%%%%%%%%%%%%%%%%%%%%%%%%%%%%%%%%%%%%%%%%
%%%%%%%%%%%%%%%%%%%%%                         %%%%%%%%%%%%%%%%%%%%%
%%%%%%%%%%%%%%%%%%%%%%%%%%%%%%%%%%%%%%%%%%%%%%%%%%%%%%%%%%%%%%%%%%%
%%%%%%%%%%%%%%%%%%%%%%%%%%%%%%%%%%%%%%%%%%%%%%%%%%%%%%%%%%%%%%%%%%%

\section {General formulation}
\setcounter{equation}{0}

\subsection {Some definitions}

The pseudoeuclidean Lie algebra $e\pq$ is the semidirect sum of the 
pseudoorthogonal Lie algebra $o(p,q)$ and an abelian algebra $T(n)$ 
of translations
\beq
e\pq=o\pq \c+ T(n), \qquad n = p+q.
\eeq

We will make use of the following matrix representation of the 
Lie algebra $e\pq$ and the corresponding Lie group $E\pq$. We introduce 
an "extended metric"
\beq 
K_{e}= \left( \begin{array}{cc} K & 0 \\
                                 0 & 0_{1} 
                \end{array} \right), \label{Ke}
\eeq
where $K$ satisfies
\bea
&K=K^{T} \in \mR^{n \times n},  \quad  n=p+q,  \quad  det K \not= 0,&  \label{K}\\
& sgn K=(p,q),  \qquad   p \geq q \geq 0. &
\eea
Here $sgnK$ denotes the signature of $K$, where $p$ and $q$ are the numbers of 
positive and negative eigenvalues, respectively.
Then $X_e \in e\pq$ and $H \in E\pq$ are represented as
\bea
& X_e(X,\alpha) \equiv X_e =\left( \begin{array}{cc} X  & \alpha^T \\
                                               0 & 0 
                              \end{array} \right),
\qquad X \in \mR^{n \times n}, \qquad \alpha \in \mR^{1 \times n} \label{H}, & \\
& H=\left( \begin{array}{cc} G & a^T \\
                            0 & 1 
           \end{array} \right),                                   
\qquad G \in \mR^{n \times n},\qquad a \in \mR^{1 \times n}, \label{G} & \\
& XK+KX^{T}=0, \quad GKG^{T}=K, \qquad X_eK_{e}+K_{e}X_e^{T}=0. \label{XK} &
\eea

The vector $\alpha \in \mR^{1 \times n}$ represents the translations.
We say that the translations are positive, negative or zero (isotropic) length if
\beq
\alpha K \alpha^T >0, \quad \alpha K \alpha^T <0, \quad \alpha K \alpha^T = 0,
\eeq
respectively.

We will be classifying maximal abelian subalgebras of the pseudoeuclidean Lie 
algebra $e(p,q)$ into conjugacy classes under the action of the pseudoeuclidean 
Lie group $E(p,q)$. Let us define some basic concepts.
\bdf
The centralizer $cent(L_{0},L)$ of a Lie algebra $L_{0}\subset L$ is a subalgebra of L 
consisting of all elements in L, commuting elementwise with $L_{0}$
\beq
cent(L_{0},L)=\{e \in L|[e,L_{0}]=0\}.
\eeq
\edf

\bdf
A maximal abelian subalgebra $L_{0}$ (MASA) of L is an abelian subalgebra, equal
to its centralizer
\beq
[L_{0},L_{0}]=0, \, \, \, \, \,  cent(L{_0},L)=L_{0}.
\eeq
\edf
\bdf
A normalizer group $N\!or(L_0,G)$ in the group $G$ of the subalgebra $L_0 \subseteq L$ is
\beq 
N\!or(L_0,G)=\{g \in G | gL_0 g^{-1} \subseteq L_0\}.
\eeq
\edf
\bdf
A splitting subalgebra $L_{0}$ of the semidirect sum
\beq
L=F {\rm \c+} N,\, \, [F,F]\subseteq F, \, \, [F,N]\subseteq N, \, \,  
[N,N]\subseteq N
\eeq
is itself a semidirect sum of a subalgebra of F and a subalgebra of N
\beq
L_{0}=F_{0} {\rm \c+} N_{0},\, \, \, F_{0}\subseteq F, \, \, \, N_{0}\subseteq N. 
\eeq
%(or conjugate to such a semidirect sum).
\edf
All other subalgebras of $L=F \c+ N$ are called 
{\it nonsplitting subalgebras}.

An {\it abelian splitting subalgebra} of $L=F\c+ N$ is a direct sum
\beq
L_{0}=F_{0}\oplus N_{0},\, \, \, \, F_{0}\subseteq F, \, \, \, \,
N_{0}\subseteq N.
\eeq 

\bdf
A maximal abelian nilpotent subalgebra (MANS) M of a Lie algebra L is a MASA, 
consisting entirely of nilpotent elements, i.e. it satisfies
\beq
[M,M]=0, \, \, \, \, \, \left[ \left[ [L,M]M\right] \ldots \right]_{m}=0
\eeq
for some finite number m (we commute M with L $m$-times). A MANS is represented by 
nilpotent matrices in any finite dimensional representation.
\edf

%%%%%%%%%%%%%%%%%%%%%%%%%%%%%%%%%%%%%%%%%%%%%%%%%%%%%%%%%%%%%%%%%%%
%%%%%%%%%%%%%%%%%%%%%%%%%%%%%%%%%%%%%%%%%%%%%%%%%%%%%%%%%%%%%%%%%%%
%%%%%%%%%%%%%%%%%%%%%                         %%%%%%%%%%%%%%%%%%%%%
%%%%%%%%%%%%%%%%%%%%%%%%%%%%%%%%%%%%%%%%%%%%%%%%%%%%%%%%%%%%%%%%%%%
%%%%%%%%%%%%%%%%%%%%%%%%%%%%%%%%%%%%%%%%%%%%%%%%%%%%%%%%%%%%%%%%%%%

\subsection {Classification strategy}

The classification of MASAs of $e(p,q)$ is based on the fact that $e(p,q)$ is 
the semidirect sum of the Lie algebra $o(p,q)$ and an abelian ideal $T(n)$ 
(the translations). We use here a procedure related to one used earlier ~\cite{kw} 
for $e(n, C)$ and \cite{conf} for $e(p,1)$. It proceeds in five steps.

\begin{zoznamnum}
\item Classify subalgebras $T(k_{+},k_{-},k_{0})$ of $T(n)$. They are 
characterized by a triplet $(k_{+},k_{-},k_{0})$, where $k_{+},k_{-}$ and $k_{0}$ 
are the number of positive length, negative length and isotropic vectors, respectively.
\item Find the centralizer $C(k_{+},k_{-},k_{0})$ of $T(k_{+},k_{-},k_{0})$ in 
$o(p,q)$
\beq C(k_{+},k_{-},k_{0})=\{X \in o(p,q)|[X,T(k_{+},k_{-},k_{0})]=0\}. \eeq
\item Construct all MASAs $M(k_{+},k_{-},k_{0})$ of $C(k_{+},k_{-},k_{0})$ and classify them 
under the 
action of normalizer $N\!or[T(k_{+},k_{-},k_{0}),G]$ of $T(k_{+},k_{-},k_{0})$
in the group $G \sim E(p,q)$.
\item Obtain a representative list of all splitting MASAs of $e(p,q)$ as  direct sums
\beq M(k_{+},k_{-},k_{0}) \oplus T(k_{+},k_{-},k_{0}) \eeq
and keep only those amongst them that are indeed maximal (and mutually inequivalent).
\item Construct all nonsplitting MASAs from splitting ones as described below in Section 5.1.
\end{zoznamnum}

%%%%%%%%%%%%%%%%%%%%%%%%%%%%%%%%%%%%%%%%%%%%%%%%%%%%%%%%%%%%%%%%%%%
%%%%%%%%%%%%%%%%%%%%%%%%%%%%%%%%%%%%%%%%%%%%%%%%%%%%%%%%%%%%%%%%%%%
%%%%%%%%%%%%%%%%%%%%%                         %%%%%%%%%%%%%%%%%%%%%
%%%%%%%%%%%%%%%%%%%%%%%%%%%%%%%%%%%%%%%%%%%%%%%%%%%%%%%%%%%%%%%%%%%
%%%%%%%%%%%%%%%%%%%%%%%%%%%%%%%%%%%%%%%%%%%%%%%%%%%%%%%%%%%%%%%%%%%

%%%%%%%%%%%%%%%%%%%%%%%%%%%%%%%%%%%%%%%%%%%%%%%%%%%%%%%%%%%%%%%%%%%
%%%%%%%%%%%%%%%%%%%%%%%%%%%%%%%%%%%%%%%%%%%%%%%%%%%%%%%%%%%%%%%%%%%
%%%%%%%%%%%%%%%%%%%%%                         %%%%%%%%%%%%%%%%%%%%%
%%%%%%%%%%%%%%%%%%%%%%%%%%%%%%%%%%%%%%%%%%%%%%%%%%%%%%%%%%%%%%%%%%%
%%%%%%%%%%%%%%%%%%%%%%%%%%%%%%%%%%%%%%%%%%%%%%%%%%%%%%%%%%%%%%%%%%%

\section{Results on MASAs of ${\mathbf o(p,q)}$}
\setcounter{equation}{0}
\subsection{General results}

Let us briefly sum up some known \cite{verop} results on MASAs of $o\pq$ that
we shall need below. We shall represent these MASAs by matrix sets $\{X,K\}$
with notations as in  (\ref{K}) \ldots (\ref{XK}).

\bdf
A MASA of $o\pq$ is called orthogonally decomposable (OD) if all matrices in the set 
 $\{X,K\}$ can be simultaneously represented by block diagonal matrices with the
same decomposition pattern. It is called orthogonally indecomposable (OID) otherwise.
\edf

\begin{proposition} \label{pmasa}
Every OD MASA of $o\pq$ can be represented by a matrix set
\bea
& X=diag(X_1, X_2, \ldots, X_k), \qquad K=diag(K_{p_1,q_1}, K_{p_2,q_2}, \ldots, K_{p_k,q_k}), &  
\nonumber \\
& X_j K_{p_j,q_j}+ K_{p_j,q_j}X_j^T=0, \qquad X_j, K_{p_j,q_j} \in \mR^{(p_j+q_j) \times (p_j+q_j)}, 
& \nonumber \\
& K_{p_j,q_j}=K_{p_j,q_j}^T, \qquad sgn K_{p_j,q_j}=(p_j,q_j),  & \\
& det K_{p_j,q_j} \not= 0, \qquad 1 \leq j \leq k, \qquad 2 \leq k \leq \left[ p+q+1 \over 2 \right], 
&  \nonumber \\
& \sum_{j=1}^{k}p_j=p, \quad \sum_{j=1}^{k}q_j=q, \quad p_1+q_1\geq p_2+q_2 \geq \ldots
 \geq p_k+q_k \geq 1, & \nonumber
\eea
where:
\begin{zoznamrom}
\item For each j, the matrix set $\{X_j,K_{p_j,q_j}\}$ represents an OID MASA of $o(p_j,q_j)$; let us call it $M_{p_j,q_j}$.
\item At most one of the MASAs  $M_{p_j,q_j}$ is a maximal abelian nilpotent subalgebra (MANS) of $o(p_j,q_j)$. In particular only one pair $(p_j,q_j)$ can satisfy $p_j+q_j=1$. The corresponding pair $\{X,K\}$ is $(0,1)$ and represents a MANS of $o(1,0)$ or $o(0,1)$.
\end{zoznamrom}
To obtain representatives of all $O(p,q)$ classes of OD MASAs of $o\pq$ we let $M_{p_j,q_j}$, for all $j$, 
run independently through all representatives of $O(p_j,q_j)$ conjugacy classes of OID MASAs of $o(p_j,q_j)$,
 subject to the restriction (ii). Conversely, each such matrix set represents a conjugacy 
class of OD MASAs of $o\pq$.
\end{proposition}

The problem of classifying MASAs of $o\pq$ is thus reduced to the classification of OID MASAs.
Under the field extension from $\mR$ to $\mC$ an OID MASA can remain OID, or become orthogonally 
decomposable. In the first case we call it {\it absolutely orthogonally decomposable} (AOID) in 
the second {\it nonabsolutely orthogonally indecomposable} (NAOID). The following types of 
orthogonally indecomposable MASAs of  $o\pq$ exist:

\noindent
1. Maximal abelian nilpotent subalgebras (MANSs). They exist for all values of $\pq$, 
$min(p,q) \geq 1$.  They are discussed below in Section 3.2. They are AOID MASAs.

\noindent
2. MASAs that are decomposable but not orthogonally decomposable (AOID but D). They stay 
OID when considered over $\mC$. They exist for
all values of $p=q\geq 1$. Their canonical form is
\beq
M=\left\{ X_{p,p}=\left( \begin{array}{cc} A & \\
                                     & -A^T
            \end{array} \right), \quad
K=\left( \begin{array}{cc}  &  I_{p} \\
                          I_{p} & 
            \end{array} \right) \right\},
\eeq
where $A= \mR I_{p} \oplus {\rm \, \, MANS \, \, of \, \,} sl(p,\mR)$. 

\noindent
3. MASAs that are indecomposable over $\mR$ but become orthogonally decomposable after field 
extension
to $\mC$  (NAOID, ID but NAID). They exist for $p=2k$, $q=2l$, $min(k,l) \geq 1$. Their canonical 
form is
$$M= \mR Q \oplus {\rm MANSs \, \, of \, \, } su(k,l), \qquad 
K=\zl \begin{array}{cc} I_{2k} &  \\  & -I_{2l} \end{array} \right),$$
\beq
Q=diag(F_2, \ldots , F_2) \in \mR^{2(k+l) \times 2(k+l)}, \qquad 
F_2=\left( \begin{array}{cc} 0 & 1 \\ -1 & 0 \end{array} \right) \label{Q}.
\eeq

\noindent
4. MASAs that are indecomposable over $\mR$ and decomposable over $\mC$ (but not orthogonally decomposable 
even over $\mC$) (OID, AOID but NAID). They exist for $p=q=2k$, $k \geq 1$. Their canonical form is
$$M= \mR Q \oplus {\rm OID \, \, but \, \, D \, \, MASAs \, \,  of \, \,} su(k,k)$$
with $Q$ as in eq.(\ref{Q}).

An exception is the case of $o(2)$, itself abelian. Thus, for $p=2$, $q=0$ or $p=0$, $q=2$, $o(2)$ is 
AOID but NAID.

\noindent
5. Decomposable MASAs that become orthogonally decomposable over $\mC$ (NAOID and D). They occur 
only for $p=q=2k$, $k \geq 1$. Their canonical form is 
\beq
M=\left\{ X= \left( \begin{array}{cc} A & \\
                                     & -A^T
            \end{array} \right), \quad
K=\left( \begin{array}{cc}  &  I_{2k} \\
                          I_{2k} & 
            \end{array} \right) \right\},
\eeq
where
$$A=\mR Q_{2k} \oplus {\rm MANSs \, \, of \, \, } sl(2k,\mC). $$

%%%%%%%%%%%%%%%%%%%%%%%%%%%%%%%%%%%%%%%%%%%%%%%%%%%%%%%%%%%%%%%%%%%
%%%%%%%%%%%%%%%%%%%%%%%%%%%%%%%%%%%%%%%%%%%%%%%%%%%%%%%%%%%%%%%%%%%
%%%%%%%%%%%%%%%%%%%%%                         %%%%%%%%%%%%%%%%%%%%%
%%%%%%%%%%%%%%%%%%%%%%%%%%%%%%%%%%%%%%%%%%%%%%%%%%%%%%%%%%%%%%%%%%%
%%%%%%%%%%%%%%%%%%%%%%%%%%%%%%%%%%%%%%%%%%%%%%%%%%%%%%%%%%%%%%%%%%%

\subsection{MANSs of $\mathbf o\pq$}

A MANS $M$ of a classical Lie algebra is characterized by its Kravchuk signature, which we will 
denote KS \cite{verop,decomp,st,krav}. It is a 
triplet of integers
\beq
(\lambda \, \, \mu \, \, \lambda), \qquad 2\lambda +\mu =n, \qquad \mu \geq 0, 
\qquad 1 \leq \lambda \leq q \leq p,
\eeq
where $\lambda$ is the dimension of the kernel of $M$, equal to the codimension of the image of $M$.
 A MANS can be transformed into 
the Kravchuk normal form
\bea
& N= \zl \begin{array}{ccc} 0 & A & Y \\
                         0 & S & -\tilde KA^T \\
                           0&0&0
         \end{array} \zr,   \qquad
K=\zl \begin{array}{ccc}  &  & I_{\lambda} \\
                          & \tilde K & \\
                          I_{\lambda} & & 
         \end{array} \zr, & \nonumber \\
& A \in \mR^{\lambda \times \mu}, \qquad Y=-Y^T \in \mR^{\lambda \times \lambda}, \qquad S\tilde K+\tilde KS^T=0, & \\
& S \in \mR^{\mu \times \mu}, \quad \tilde K=\tilde K^T \in \mR^{\mu \times \mu}, \quad
sgn\tilde K=(p-\lambda,q-\lambda) & \nonumber
\eea
and S nilpotent.

There are two types of MANS of $o\pq$:
\begin{zoznamrom}
\item
Free-rowed MANS. The first row of $A$  has $\mu$ free real entries. 
All other 
entries in $A$ and $S$ depend linearly on those $\mu$ free entries.
\item Non-free-rowed MANS. Any combination of rows of $A$ contains less than $\mu$ free real entries.
\end{zoznamrom}

The results on free-rowed MANS of $o\pq$  \cite{verop} are stated in the following proposition.
\begin{proposition} \label{pmans}
A representative list of $O\pq$ conjugacy classes of free-rowed MANSs of $o\pq$ with Kravchuk 
signature {\it ($\lambda$ $\mu$ $\lambda$)} is given by the matrix sets
\bea
N= \zl \begin{array}{ccc}0 & A & Y \\
                         0 & 0 & -\tilde KA^T \\
                           0&0&0
         \end{array} \zr,   \qquad
K=\zl \begin{array}{ccc}  &  & I_{\lambda} \\
                          & \tilde K & \\
                          I_{\lambda} & & 
         \end{array} \zr,  \label{manso} \\
A= \zl \begin{array}{c} \alpha Q_1 \\ \alpha Q_2 \\ \vdots \\ \alpha Q_{\lambda}
          \end{array} \zr, \qquad 
\begin{array}{c} \alpha \in \mR^{1 \times \mu} \\
                  Y=-Y^T \in \mR^{\lambda \times \lambda}, \label{AY}
\end{array} \\
Q_i \in \mR^{\mu \times \mu}, \quad Q_i\tilde K=\tilde KQ_i^T, \quad [Q_i,Q_j]=0, \label{Qi} \\  
Q_1=I, \quad 
 TrQ_i=0, \qquad 2 \leq i \leq \lambda.  \nonumber
\eea 
The entries in $\alpha$ and Y are free. The matrices $Q_i$ are fixed and form an abelian subalgebra of the 
Jordan algebra $jo(p-\lambda,q-\lambda)$. In the case $\lambda =2$ we must have $Q_2 \not= 0$. There exists a
$\lambda_{1} \in \mZ, 1 \leq \lambda_{1} \leq \lambda$ such that $Q_1, \ldots, Q_{\lambda_{1}}$ are
linearly independent and $Q_{\nu} =0$, $\lambda_{1}+1 \leq \nu \leq \lambda$.
\end{proposition}
Proofs of the Propositions \ref{pmasa} and \ref{pmans} and details about MASAs of $o\pq$ 
are given in Ref.~\cite{verop}. The results on non-free-rowed MANS of $o\pq$ are less complete and we shall not reproduce them here \cite{verop}.

%%%%%%%%%%%%%%%%%%%%%%%%%%%%%%%%%%%%%%%%%%%%%%%%%%%%%%%%%%%%%%%%%%%
%%%%%%%%%%%%%%%%%%%%%%%%%%%%%%%%%%%%%%%%%%%%%%%%%%%%%%%%%%%%%%%%%%%
%%%%%%%%%%%%%%%%%%%%%                         %%%%%%%%%%%%%%%%%%%%%
%%%%%%%%%%%%%%%%%%%%%%%%%%%%%%%%%%%%%%%%%%%%%%%%%%%%%%%%%%%%%%%%%%%
%%%%%%%%%%%%%%%%%%%%%%%%%%%%%%%%%%%%%%%%%%%%%%%%%%%%%%%%%%%%%%%%%%%

%%%%%%%%%%%%%%%%%%%%%%%%%%%%%%%%%%%%%%%%%%%%%%%%%%%%%%%%%%%%%%%%%%%
%%%%%%%%%%%%%%%%%%%%%%%%%%%%%%%%%%%%%%%%%%%%%%%%%%%%%%%%%%%%%%%%%%%
%%%%%%%%%%%%%%%%%%%%%                         %%%%%%%%%%%%%%%%%%%%%
%%%%%%%%%%%%%%%%%%%%%%%%%%%%%%%%%%%%%%%%%%%%%%%%%%%%%%%%%%%%%%%%%%%
%%%%%%%%%%%%%%%%%%%%%%%%%%%%%%%%%%%%%%%%%%%%%%%%%%%%%%%%%%%%%%%%%%%

\section{Splitting MASAs of ${\mathbf e(p,q)}$}
\setcounter{equation}{0}
\subsection{General comments on MASAs of ${\mathbf e\pq}$}

A MASA of  $e\pq$ will be represented by a matrix set $\{X_e,K_e\}$
\bea 
& X_e=\left( \begin{array}{ccccccc} N & & & & & & \xi^T \\
                                      & X_{p_1,q_1} & & & & & \delta_1^T \\
                                      & & \ddots & & & & \vdots \\
                                      & & & X_{p_j,q_j} & & &  \delta_j^T \\
                                     & & & & 0_{\kp} & & x^T \\
                                     & &  & & & 0_{\km} & y^T \\
                                       & & & & &  & 0_1
            \end{array} \right), \label{X} & \\
& K_e=\left( \begin{array}{ccccccc}  K_0 & & & & & & \\
                                     & K_{p_1,q_1} & & & & & \\
                                  & & \ddots & & & &  \\
                                   & & & K_{p_j,q_j} & & & \\
                                  & & & & I_{\kp} & &  \\
                                  & & & & & -I_{\km} &  \\
                                  & & & & & & 0_1
         \end{array} \right), & \label{Ke1} \\
&  p=p_0+\kz + \sum_{i=1}^{j}p_i + \kp,  \qquad   q=q_0+ \kz + \sum_{i=1}^{j}q_i + \km, &
\eea
where $M_{p_i,q_i}=\{X_{p_i,q_i}, K_{p_i,q_i}\}$, $i=1, \ldots j$ is an OID MASA of 
$o(p_i,q_i)$, 
that is not a MANS. The vector $\xi$ has the following form
\bea
\xi =\left( \begin{array}{c} z^T \\ \beta^T \\ \gamma^T \end{array} \right), \qquad
\begin{array}{l} z, \gamma \in \mR^{1 \times \kz} \\ \beta \in \mR^{1 \times (p_0+q_0)} 
\end{array} \label{xit}
\eea
and $N$ is a MANS of $o(p_0+\kz,q_0+\kz)$ with Kravchuk signature $ (\kz \, \, p_0\! \!+\! \!q_0 \, \, \kz)$ 
and is given by
\bea 
& N=\left( \begin{array}{ccc} 0_{\kz} & A & Y \\
                             0 &     S & -K_{p_0,q_0}A^T \\
                             0 & 0 &   0_{\kz}
         \end{array} \right), 
K_0=\left( \begin{array}{ccc} 0 & 0 & I_{\kz} \\
                             0 &   K_{p_0,q_0} & 0 \\
                              I_{\kz} & 0 &   0
         \end{array} \right) & \label{mans} \\ 
&  Y=-Y^T, \qquad SK_{p_0,q_0}+K_{p_0,q_0}S^T=0 & \nonumber \\
&A \in \mR^{\kz \times (p_0+q_0)}, \qquad 
S \in \mR^{ (p_0+q_0) \times (p_0+q_0)}, \qquad Y \in \mR^{\kz \times \kz},& \\
&K_{p_0,q_0}=K^T_{p_0,q_0}, \qquad sgnK_{p_0,q_0}=(p_0,q_0)& \nonumber
\eea

The entries in $z, x$ and $y$ are free and represent the positive, negative and zero 
length translations contained in $T(\kp, \km, \kz)$. The entries in $\beta,\gamma$ 
and $\delta_i$ are linearly dependent on the free entries in $A,Y$ and $X_{p_i,q_i}$. 
If they are nonzero (and cannot be annulled by an $E\pq$ transformation), we have 
a nonsplitting MASA. This case will be discussed in Section 5.

%%%%%%%%%%%%%%%%%%%%%%%%%%%%%%%%%%%%%%%%%%%%%%%%%%%%%%%%%%%%%%%%%%%
%%%%%%%%%%%%%%%%%%%%%%%%%%%%%%%%%%%%%%%%%%%%%%%%%%%%%%%%%%%%%%%%%%%
%%%%%%%%%%%%%%%%%%%%%                         %%%%%%%%%%%%%%%%%%%%%
%%%%%%%%%%%%%%%%%%%%%%%%%%%%%%%%%%%%%%%%%%%%%%%%%%%%%%%%%%%%%%%%%%%
%%%%%%%%%%%%%%%%%%%%%%%%%%%%%%%%%%%%%%%%%%%%%%%%%%%%%%%%%%%%%%%%%%%

\subsection{Basic results on splitting MASAs}
In this section we shall construct all splitting MASAs of $e\pq$. 
\bth
Every splitting MASA of  $e\pq$ is characterized by a partition 
\beq
\begin{array}{ccccc} \label{cond}
& p=p_0+\kp +\kz + \sum_{i=1}^{j}p_i,& \qquad  & q=q_0+\km +\kz + \sum_{i=1}^{j}q_i &  \\
& \kz +\kp+\km \not= p+q-1,& \qquad & 0 \leq \kz \leq q.  & 
\end{array}
\eeq
A representative list of $E\pq$ conjugacy classes of MASAs of $e\pq$ is given 
by the matrix sets $\{X_e,K_e\}$ of eq.(\ref{X}) and (\ref{Ke1}) with  
\bea 
\delta_i=0, \quad i=1, \ldots j, \qquad
\xi=\left( \begin{array}{c} z^T \\ 0 \\ 0 \end{array} \right). \label{splitt}
\eea
If $\kz = 0$ then the MANS $N$ is absent.  $M_{p_i,q_i}$ is an orthogonally 
indecomposable MASA of $o(p_i,q_i)$ which is not a MANS. Running through all 
possible partitions, all MANSs $\{N,K_0\}$ and all MASAs $M_{p_i,q_i}$ we obtain 
a representative list of all splitting MASAs of $e\pq$.
\enth

{\it Proof}:
We start by choosing a subalgebra $T(\kp,\km,\kz)$. Calculating the centralizer of 
$T(\kp,\km,\kz)$ in 
$o\pq$ gives us
\bea
& C(\kp,\km,\kz)=\left( \begin{array}{ccc} \tilde M & &  \\
                                            & 0_{\kp} & \\
                                             & & 0_{\km}
            \end{array} \right), \qquad
K=\left( \begin{array}{ccc} \tilde K & &  \\
                            & I_{\kp} & \\
                           & & -I_{\km}
            \end{array} \right), & \\
&  sgn \tilde K=(p-\kp,q-\km).&  \nonumber
\eea
$\tilde M$ is a subalgebra of $o(p-\kp,q-\km)$ which commutes with the translations 
corresponding to $\xi = (z,0)$, $\xi \in \mR^{1 \times
(p+q-\kp-\km)}$, 
$z \in \mR^{1 \times \kz}$, and  with no other translations.
To obtain a MASA of $e\pq$ we must complement $T(\kp,\km,\kz)$ by a MASA $F(\kp,\km,\kz)$ 
of the centralizer $C(\kp,\km,\kz)$. $F(\kp,\km,\kz)$ must not commute with any further 
translations, 
 hence $F(\kp,\km,\kz)$ is either a MANS of $o(p-\kp,q-\km)$ with KS
($\kz$, \, $p\!-\!\kp\!-\!\kz\!+\!q\!-\!\km\!-\!\kz$, \, $\kz$) or an orthogonally decomposable MASA 
containing a MANS $N$
with KS ($\kz$ $\mu$ $\kz$). For 
$\kz =0$ the MANS $N$ is absent. This leads to eq.~({\ref{splitt}}) and each 
$M_{p_i,q_i}=\{X_{p_i,q_i},K_{p_i,q_i}\}$ is an OID MASA of $o(p_i,q_i)$ of 
the type 2,3,4, or 5, listed in Section 3.1.
$\hfill\Box$

%%%%%%%%%%%%%%%%%%%%%%%%%%%%%%%%%%%%%%%%%%%%%%%%%%%%%%%%%%%%%%%%%%%
%%%%%%%%%%%%%%%%%%%%%%%%%%%%%%%%%%%%%%%%%%%%%%%%%%%%%%%%%%%%%%%%%%%
%%%%%%%%%%%%%%%%%%%%%                         %%%%%%%%%%%%%%%%%%%%%
%%%%%%%%%%%%%%%%%%%%%%%%%%%%%%%%%%%%%%%%%%%%%%%%%%%%%%%%%%%%%%%%%%%
%%%%%%%%%%%%%%%%%%%%%%%%%%%%%%%%%%%%%%%%%%%%%%%%%%%%%%%%%%%%%%%%%%%

%%%%%%%%%%%%%%%%%%%%%%%%%%%%%%%%%%%%%%%%%%%%%%%%%%%%%%%%%%%%%%%%%%%
%%%%%%%%%%%%%%%%%%%%%%%%%%%%%%%%%%%%%%%%%%%%%%%%%%%%%%%%%%%%%%%%%%%
%%%%%%%%%%%%%%%%%%%%%                         %%%%%%%%%%%%%%%%%%%%%
%%%%%%%%%%%%%%%%%%%%%%%%%%%%%%%%%%%%%%%%%%%%%%%%%%%%%%%%%%%%%%%%%%%
%%%%%%%%%%%%%%%%%%%%%%%%%%%%%%%%%%%%%%%%%%%%%%%%%%%%%%%%%%%%%%%%%%%

\section{Nonsplitting MASAs of ${\mathbf e\pq}$}
\setcounter{equation}{0}

\subsection{General comments}
First we describe the general procedure for finding nonsplitting MASAs of $e(p,q)$.

Every nonsplitting MASA $M(\kp, \km, \kz)$ of $e(p,q)$ is obtained from a splitting 
one by the following procedure: 

\begin{zoznamnum}
\item Choose a basis for $F(\kp, \km, \kz)$ and $T(\kp, \km, \kz)$ e.g. 
$F(\kp, \km, \kz) \sim \\
\{B_{1}, \ldots , B_{J}\}$, $T(\kp, \km, \kz) \sim \{X_{1}, 
\ldots X_{L}\}$.
\item  Complement the basis of $T(\kp, \km, \kz)$ to a basis of $T(n)$.
\begin{eqnarray*}
T(n) / T(\kp, \km, \kz) = \{Y_{1}, \ldots , Y_{N}\}, \qquad L+N=n.
\end{eqnarray*}
\item Form the elements 
\beq
\tilde B_{a}= B_{a} + \sum_{j=1}^{N} \tilde \alpha_{aj} Y_j, \qquad a=1, \ldots, J,
\eeq
where the constants $\tilde \alpha_{aj}$ are such that $\tilde B_{a}$ form an 
abelian Lie algebra
$[\tilde B_{a},\tilde B_{b}]=0$. This provides a set of linear equations for the 
coefficients $\tilde \alpha_{aj}$. Solutions  $\tilde \alpha_{aj}$ are called 
1-cocycles and they provide  abelian subalgebras $\tilde M (\kp, \km, \kz) \sim 
\{\tilde B_{a}, X_b\} \subset e(p,q)$.
\item Classify the subalgebras $\tilde M (\kp, \km, \kz)$ into conjugacy classes 
under the action of the group $E(p,q)$. This can be done in two steps.
\begin{zoznamrom}
\item Generate trivial cocycles $t_{aj}$, called coboundaries, using the translation group $T(n)$
\beq
e^{\theta_jP_j}\tilde B_{a}e^{-\theta_jP_j}=\tilde B_{a}+\theta_j[P_j,\tilde B_{a}]=\tilde B_{a}+\sum_jt_{aj}P_j.
\eeq
The coboundaries should be removed from the set of the cocycles. If we have $\tilde 
\alpha_{aj} =t_{aj}$ for all $(a,j)$ the algebra is splitting ({\it i.e.} equivalent to a splitting one).
\item Use the normalizer of the original splitting subalgebra in the group $O(p,q)$ to further 
simplify and classify the nontrivial cocycles.

\end{zoznamrom}
\end{zoznamnum}

The general form of a nonsplitting MASA of $e\pq$ is $M_e=\{X_e,K_e\}$ given by eq. 
(\ref{X}) and (\ref{Ke1}). Requiring commutativity $[X_e,X'_e]=0$ leads to
\beq
\begin{array}{ccccc}
&X_{p_i,q_i}\delta^{'T}_i &=& X'_{p_i,q_i}\delta_i^T  &  \\
& N\xi'^T&=&N'\xi^T. & 
\end{array} \label{**}
\eeq

From the eq.(\ref{**}) we see that the entries in $\delta_i$ depend linearly only on 
$X_{p_i,q_i}$, {\it i.e.} only on the MASA  $M_{p_i,q_i}$ of $o(p_i,q_i)$.

Each $M_{p_i,q_i}$ belongs  to one of the four types of OID MASAs of $o(p_i,q_i)$ which 
were listed in the 
Section 3.1 -  AOID but D MASAs, AOID but NAID MASAs, NAOID ID but NAID MASAs or 
NAOID but D MASAs.

We will make  use of the following result:
\bl
If M is a MASA of $o(p,q)$ when considered over $\mR$, then it will also be a 
MASA of $o(n, {\mC})$, $n=p+q$, when considered over ${\mC}$.
\el

If any of the vectors $\delta_i$ were  non zero  then after field extension we 
would obtain a nonsplitting MASA of $e(n,{\mC})$ of a type that does not exist \cite{kw}.
This implies that all of the $\delta_i's$ are zero.

Any further study of nonsplitting MASAs of $e\pq$ is reduced to studying the matrices
\beq 
X_e=\left( \begin{array}{ccccccc} N & & & & &  & \xi^T\\
                                    & M_{p_1,q_1} & & & & & 0 \\
                                      & & \ddots & & & & \vdots \\
                                       & & & M_{p_j,q_j} & & &  0 \\
                                      & & & & 0_{\kp} & & 0 \\
                                      & &  & & & 0_{\km} & 0 \\
                                       & & & & &  & 0_1
            \end{array} \right) \label{delta0}
\eeq
with $\xi$ and $N$ as in eq.(\ref{xit}) and (\ref{mans}), respectively.
Further, we can see from eq.(\ref{**}) and (\ref{delta0}) that the study of nonsplitting 
MASAs is in fact reduced to 
the study of nonsplitting MASAs of $e(p_0+\kz,q_0+\kz)$ for which the projection onto the subalgebra 
$o(p_0+\kz,q_0+\kz)$ is a MANS with Kravchuk signature $(\kz \, \, \mu \, \, \kz)$, $\mu=p_0+q_0$. 
Further classification is performed under the group $E(p_0+\kz,q_0+\kz)$.

The MASAs of $e(p_0+\kz,q_0+\kz)$  to be considered will thus be represented by the matrix sets $\{X_e,K_e\}$
\beq
X_e= \left( \begin{array}{cccc} 0_{\kz} & A & Y & z^T \\
                                  0 & S & -K_{p_0,q_0}A^T & \beta^T \\
                                  0 & 0 & 0_{\kz} & \gamma^T \\
                                 0 & 0 & 0 & 0
                   \end{array} \right) ,  \qquad 
K_e=\left( \begin{array}{cccc}  &  & I_{\kz} &  \\
                                  & K_{p_0,q_0} &  &  \\
                                  I_{\kz} &  & & \\
                                  &  &  & 0_1
             \end{array} \right), \label{Zt}
\eeq 
where $Y=-Y^T$, and $\beta \in \mR^{1 \times \mu}$, $\gamma \in \mR^{1 \times \kz}$ depend 
linearly on the free entries in $A$ and $Y$. Using the commutativity $[X_e,X'_e]=0$ we obtain
\bea
A\beta'^T+Y\gamma'^T & = & A'\beta^T + Y'\gamma^T \label{*}\\
S\beta'^T-K_{p_0,q_0}A^T\gamma'^T & = & S'\beta^T-K_{p_0,q_0}A'^T\gamma^T \nonumber
\eea 

The translations 
\beq
\Pi= \left( \begin{array}{cccc} 0_{\kz} & 0 & 0 & 0 \\
                                 0 &  0_{p_0,q_0} & 0 & \tau^T \\
                                  0 & 0 & 0_{\kz} & \zeta^T \\
                                 0 & 0 & 0 & 0_1
             \end{array} \right), \qquad
\tau \in \mR^{1 \times \mu}, \zeta \in \mR^{1 \times \kz}. \label{Pit}
\eeq
will be used to remove coboundaries from $\beta$ and $\gamma$ and the remaining
cocycles will be classified under the action of the normalizer of the MANS $N$ in the
group $O(p_0+\kz,q_0+\kz)$.

The situation will be very different for free-rowed and non-free-rowed MANS of $o(p_0+\kz,q_0+\kz)$.
The two cases will be treated separately.

%%%%%%%%%%%%%%%%%%%%%%%%%%%%%%%%%%%%%%%%%%%%%%%%%%%%%%%%%%%%%%%%%%%
%%%%%%%%%%%%%%%%%%%%%%%%%%%%%%%%%%%%%%%%%%%%%%%%%%%%%%%%%%%%%%%%%%%
%%%%%%%%%%%%%%%%%%%%%                         %%%%%%%%%%%%%%%%%%%%%
%%%%%%%%%%%%%%%%%%%%%%%%%%%%%%%%%%%%%%%%%%%%%%%%%%%%%%%%%%%%%%%%%%%
%%%%%%%%%%%%%%%%%%%%%%%%%%%%%%%%%%%%%%%%%%%%%%%%%%%%%%%%%%%%%%%%%%%
\subsection{Nonsplitting MASAs of ${\mathbf e(p_0+\kz,q_0+\kz)}$ related to 
free-rowed MANSs}
\def\zl{\left(}
\def\zr{\right)}

Let $N$ be a free-rowed MANS  of $o(p_0+\kz,q_0+\kz)$. The corresponding 
nonsplitting MASAs of 
$e(p_0+\kz,q_0+\kz)$ can be represented as follows.

\bth \label{tfr}
A nonsplitting MASA of $e\pq$ must contain a MANS of $o(p_0+\kz,q_0+\kz)$ with 
$1 \leq \kz \leq q$, $min(p_0+\kz,q_0+\kz) \geq 1$.
All nonsplitting MASAs of $e(p_0+\kz,q_0+\kz)$ for which the projection onto 
$o(p_0+\kz,q_0+\kz)$
is a free rowed MANS $N$ with Kravchuk signature  $(\kz \, \, \mu \, \, \kz)$, 
$\mu=p_0+q_0$ can 
be represented by the matrix sets $\{X_e,K_e\}$ of eq.(\ref{Zt}) with $S=0$ 
and $A$ and $Y$ 
as in eq.(\ref{AY}). \\
1. For $\kz \geq 3$ we have 
\beq
\beta= a \Lambda, \qquad \gamma=0
\eeq
$\Lambda \in \mR^{\mu \times \mu}$ satisfies the following conditions:
\beq 
\Lambda = \Lambda^T,  \qquad \qquad 
Q_j\Lambda K_{p_0,q_0}^{-1}=\Lambda K_{p_0,q_0}^{-1} Q_j. \label{lQ}
\eeq
2. $\kz=2$, $\mu \geq 2$.
$\Lambda$ satisfies eq.(\ref{lQ}) for $j=2$ and
\bea
\beta = \alpha\Lambda + y \rho, \qquad \gamma= \left( \begin{array}{c} 0 \\ \alpha \rho^T \end{array} \right), 
\qquad \rho=(1, 0, \ldots , 0) 
\eea
for Q following
\bea
Q=\left( \begin{array} {cccc} 0 & 1 & & \\
                               &  0 & & \\
                               & & \ddots & \\
                                & & & 0
            \end{array} \right), \qquad
K_{p_0,q_0}= \left( \begin{array} {ccc} 0 & 1 &  \\
                               1 &  0 &  \\
                               & &  K_{p_0-1,q_0-1}
            \end{array} \right). \label{QQ}
\eea
For all the other $Q$ 
\bea
\beta= \alpha \Lambda, \qquad \gamma= \left( \begin{array}{c} 0 \\ 0 \end{array} \right).
\eea
3. $\kz=2$, $\mu=1$
\beq
\beta=y \rho, \qquad \gamma^T= \left( \begin{array}{c} 0 \\ a \rho+ p_2 y \end{array} \right),
\eeq
where $(\rho,p_2)$ is $(1,0), (0,1),$ or $(1,1)$. \\
4. $\kz=2$, $\mu=0$, there is no $\beta$ and we have
\beq
\gamma^T= \left( \begin{array}{c}  y \\ 0 \end{array} \right). \label{ss}
\eeq
5. $\kz=1$, $\mu \geq 2$
\beq
\beta = \alpha \Lambda, \qquad \Lambda = \Lambda ^T, \qquad \gamma=0.
\eeq
6. $\kz=1$, $\mu=1$
\beq
\beta =0, \qquad \gamma =a.
\eeq
The case $\kz=1$, $\mu=0$ is not allowed. \\
Two free-rowed non-splitting MASAs of $e(p_0+\kz,q_0+\kz)$, $M(p_0,q_0,\kz,\Lambda)$ and
$M'(p_0,q_0,\kz,\Lambda')$, are $E(p_0+\kz,q_0+\kz)$ conjugated (for cases 1 and 5) if 
the matrices $\Lambda$, $\Lambda'$ characterizing them satisfy:
\beq
\Lambda'={1 \over g_1} G_2 ( \Lambda - \sum_{k=1}^{\kz} \theta_{k} Q_k K_{p_0,q_0}) G_2^T \label{lambda}
\eeq 
for some $g_1,g_j \in \mR$, $\theta_k \in \mR$, $G_2 \in o(p_0,q_0)$ such that
\beq
Q_j={1 \over g_1} g_j G_2 Q_j G_2^{-1} \label{qj}.
\eeq

\enth

{\it Proof:}
1. $\kz \geq 3$
We start from  a free-rowed MANS in eq.(\ref{Zt}). Requiring commutativity $[X_e, X'_e]=0$ 
leads to the following equations
\beq
\begin{array}{rcl}
(\alpha Q_j) \beta^{'T} + y_{ja}\gamma'_a&=&(\alpha' Q_j) \beta^{T} + y'_{ja}\gamma_a  \\
(Q_j \alpha^T)\gamma'_j&=&(Q_j \alpha'^T)\gamma_j. \label{a}
\end{array}
\eeq
The entries in $\beta, \gamma$ are linearly dependent on those in $Y$ and $\alpha$, {\it i.e.}
\beq
\begin{array}{cclccc}
\beta &= & \alpha \Lambda + \sum_{1 \leq i < k \leq \kz} y_{ik}\rho_{ik}, & & &
\Lambda \in \mR^{\mu \times \mu},  \rho_{ik} \in \mR^{1 \times \mu} \\
\gamma&=& \alpha W + \sum_{1 \leq i < k \leq \kz} y_{ik}P_{ik}, & & &
W \in \mR^{\mu \times \kz},  P_{ik} \in \mR^{1 \times \kz} 
\end{array}
\eeq
We substitute $\beta$ and $\gamma$ into eq.~(\ref{a}) and 
compare coefficients of $\alpha_i \alpha'_j$, for $i$ and $j$ fixed. First consider the case $j=1$.
We obtain that
\bea
& \Lambda = \Lambda^T; \qquad P_{ik,a}=0, \quad 2\leq i < k , \quad 1 < a; \qquad P_{1k,a}=P_{1a,k},& 
\nonumber \\
& \rho_{ik}=0, \quad 2 \leq i <k; \qquad W_{a}=\rho_{1a}, \quad a \geq 2, & \\
& Q_j\Lambda K_{p_0,q_0}^{-1}=\Lambda K_{p_0,q_0}^{-1} Q_j. & \nonumber
\eea
For $j=2$ we obtain
\beq 
\begin{array}{rclcrcl}
P_{ik,1}& =& 0 & \quad 3 \leq i<k, & \qquad P_{12,a}& =& -P_{2a,1} \\
\rho_{1k}& =& 0 & \quad k \geq 3, & \qquad W_1&=&-Q_2\rho_{12}^T.
\end{array} 
\eeq
And for $j=3$ we get
\beq
W=0, \qquad \rho_{ik}=0, \qquad P_{ik}=0 \qquad {\rm for} \quad  \kz \geq 3,
\eeq
Using the translations we obtain the coboundaries $\theta_i$
\bea
e^{\theta_i P_i} Z e^{- \theta_i P_i} = Z - \theta_i [Z, P_i].
\eea
This leads to replacing $\Lambda$ by
\bea
\Lambda ' = \Lambda - \sum_{k=1}^{\kz} \theta_k Q_k K_{p_0,q_0}.
\eea
All $\theta_i$ are free and can be used to remove all coboundaries. In particular if  $K_{p_0,q_0}$ 
is chosen to satisfy $TrK_{p_0,q_0} \neq 0$ we can use $\theta_1$ to make $\Lambda$ traceless. 
Equation (\ref{lambda}) corresponds to transformations of  $\Lambda$ using the normalizer of $N$ in $E\pq$.

\noindent 2. $\kz=2$, $\mu \geq 2$ \\
Here there is only one matrix $Q=Q_2$, the vector $\gamma$ is  $\gamma =(\gamma_1,\gamma_2)$ 
and $Y=\zl \begin{array}{cc} 0 & y \\ -y &0 \end{array} \zr$. We have
\bea
& \beta = \alpha \Lambda +y\rho, \qquad \rho \in \mR^{1 \times \mu} & \\
& \gamma_1=\alpha w_1^T+p_1y, \quad \gamma_2=\alpha w_2^T+p_2y, \quad w_1,w_2 \in \mR^{1 \times \mu}, 
\quad p_1,p_2 \in \mR &
\eea
From the  $[X_e, X'_e]=0$ we obtain that
\bea
&\Lambda=\Lambda^T, \qquad  Q\Lambda K_{p_0,q_0}^{-1}=\Lambda K_{p_0,q_0}^{-1} Q &\\
&\beta = \alpha\Lambda + y \rho, \qquad  \gamma= \left( \begin{array}{c} -\alpha Q \rho^T \\ \alpha \rho^T \end{array} \right). &
\eea
Equation (\ref{a}) for $j=2$ leads to
\bea
[Q^T(\alpha^T\alpha'-\alpha'^T\alpha)+(\alpha'^T\alpha-\alpha^T\alpha')Q]\rho^T=0. \label{qqq}
\eea
Writing eq.(\ref{qqq}) in components and choosing $\alpha$ and $\alpha'$ such that $\alpha_a=1, \alpha'_b=1$ and all other components vanish, we obtain
\beq
(Q^T)_{ia}\rho_b-(Q^T)_{ib}\rho_a+\sum_{k=1}^{\mu}(\delta_{ib}Q_{ak}-\delta_{ia}Q_{bk})\rho_k=0, \qquad \forall i,a,b. \label{ind}
\eeq
This provides us with two types of relations
\bea
Q_{ai}\rho_b-Q_{bi}\rho_a=0, \qquad a \neq i, \quad b\neq i \label{qa} \\
-Q_{ii}\rho_a+Q_{ai}\rho_i+\sum_{k=1}^{\mu} Q_{ak}\rho_k=0, \qquad a\neq i. \label{qb}
\eea
The matrix $Q$ is block diagonal,
\bea
\begin{array}{c}Q=diag(J_1,J_2, \ldots, J_r), \qquad \sum_{i=1}^{r} dimJ_i =\mu \\  \\
dimJ_1 \geq dimJ_2 \geq \ldots \geq dimJ_r \geq 1, \end{array}
\eea
where each $J_i$ is an indecomposable element of a Jordan algebra $jo(p_i,q_i)$,
 $p_i+q_i=dimJ_i$ (see e.g Ref. \cite{dj}). The matrix $K_{p_0,q_0}$ has the same block structure.
 Possible forms of elementary blocks in $Q$ are
\beq
\begin{array}{c} J_i(q_i)= \zl \begin{array}{cccc} q_i & 1 & & \\
                                  & q_i & 1 & \\
                                   & & \ddots & 1 \\
                                   & & & q_i
\end{array} \zr, \\
J_i(r_i+s_i)= \zl \begin{array}{ccccccc} r_i & s_i & 1 & 0 & & & \\
                                         -s_i & r_i & 0 & 1 & & & \\
                                              & & \ddots & & \ddots & & \\
                                          & & & & & 1 & 0 \\
                                          & & & & & 0 & 1 \\
                                          & & & & & r_i & s_i \\
                                          & & & & & -s_i & r_i 
               \end{array} \zr. \end{array}
\eeq
After complexification the second type of block reduces to the first one, so it actually 
suffices to consider the first type of block only (see Lemma 5.1).

Let us first assume $dimJ_1 \geq3$. Writing relation (\ref{qb}) for $i=1$ and 
$2 \leq a\leq r$ 
we obtain $\rho_3=\rho_4=\ldots =\rho_{\mu} =0$. Taking $a=1, i=2$ in (\ref{qa}) we then 
obtain $\rho_2=0$. Taking $a=1, b=2, i=3$ in (\ref{qb}) we obtain $\rho_1=0$. Thus, if the 
largest block $J_1(q)$ satisfies $dimJ_1(q) \geq 3$, we have $\rho=0$.

Now let us assume $dimJ_1(q)=2$ so that all other blocks have dimension 2 or 1. By the same 
argument we have $\rho_3=\rho_4=\ldots =\rho_{\mu} =0$ and also $\rho_2 =0$. If $Q$ has the form
(\ref{QQ}), then all relations (\ref{qa}) and (\ref{qb}) are satisfied and $\rho_1$ remains
free. If any of the other diagonal elements, say $Q_{33}$ is not zero, then relation 
(\ref{qb}) for $i=3, a=1$ implies $\rho_1 =0$. If we have $q \neq 0$ in $J_1(q)$, then 
at least one other diagonal element of $Q$ must satisfy $Q_{aa} \neq 0, a \geq 3$, since 
we have $Tr Q=0$. 

Finally, let $Q$ be diagonal. We have $Q \neq 0, TrQ=0$, hence at least two diagonal elements 
are nonzero. Relations (\ref{qa}) and (\ref{qb}) then imply $\rho_i=0, i=1, \ldots \mu$.

Using the normalizer $G=diag(g_1,g_2,G_2,g_1^{-1},g_2^{-1})$ we normalize $\rho_1$ to $\rho_1=1$ 
for $\rho_1 \neq 0$.

\noindent 3.  $\kz=2$, $\mu= 1$ \\
There is no matrix $Q$ and we have
\bea
&\beta=\lambda a+\rho y, \qquad \lambda \in \mR & \\
& \gamma_1=aw_1+p_1y, \qquad \gamma_2=aw_2+p_2y \qquad
w_1,w_2,p_1,p_2 \in \mR. \nonumber
\eea
Condition  $[X_e, X'_e]=0$ implies $w_1=0, p_1=0$
and after removing the coboundaries we obtain
\bea
\beta =\rho y, \qquad \gamma_1=0, \qquad \gamma_2=a \rho+p_2 y.
\eea
Using the normalizer $G=diag(g_1, g_2, g_3, g_4, g_5,1)$, satisfying $G\tilde K_0 G^T = \tilde K_0$, 
we can normalize $(\rho, p_2)$ to one of the following: $(1,0), (1,1), (0,1)$.

\noindent 4. $\kz=2$, $\mu= 0$ \\
Using the normalizer $G=diag(g_1, G_2, {1 \over g_1},1)$ we  obtain eq.~(\ref{ss}). 

\noindent 5. $\kz=1$, $\mu \geq 2$ \\
In this case $Y=0$ and $A=\alpha \in \mR^{1 \times \mu}$ in eq.~(\ref{Zt}). Then we have
\beq
\beta = \alpha \Lambda , \qquad \gamma=\alpha w^T, \quad \beta \in \mR^{1 \times \mu}, \quad \gamma \in \mR.
\eeq
From the  $[X_0, X'_0]=0$ we obtain that
\beq
\Lambda=\Lambda^T, \qquad w=0. 
\eeq
Removing the coboundaries leads to replacing $\Lambda$ by
\bea
\Lambda ' = \Lambda - \theta K_{p_0,q_0},
\eea
where $\theta$ can be chosen to annul trace of $\Lambda$ (if $TrK_{p_0,q_0} \neq 0$).

\noindent 6. $\kz=1$, $\mu=1$ \\
The proof is trivial and can be found in Ref.\cite{conf}. 

Using the normalizer of the splitting MASA (\ref{splitt})in the group $E(p_0+\kz,q_0+\kz)$
we can simplify $\Lambda$ further. The normalizer is represented by block diagonal matrices
\beq
G=diag(G_1,G_2,G_1^{-1},1). 
\eeq
Choosing $G_1=diag(g_1, \ldots, g_{\kz})$, $G_2$ satisfying $G_2  K_{p_0,q_0} G_2^T= 
K_{p_0,q_0}$ 
leads to equations (\ref{lambda}) and (\ref{qj}).

This completes the proof of the Theorem 5.1.
$\hfill\Box$

%%%%%%%%%%%%%%%%%%%%%%%%%%%%%%%%%%%%%%%%%%%%%%%%%%%%%%%%%%%%%%%%%%%
%%%%%%%%%%%%%%%%%%%%%%%%%%%%%%%%%%%%%%%%%%%%%%%%%%%%%%%%%%%%%%%%%%%
%%%%%%%%%%%%%%%%%%%%%                         %%%%%%%%%%%%%%%%%%%%%
%%%%%%%%%%%%%%%%%%%%%%%%%%%%%%%%%%%%%%%%%%%%%%%%%%%%%%%%%%%%%%%%%%%
%%%%%%%%%%%%%%%%%%%%%%%%%%%%%%%%%%%%%%%%%%%%%%%%%%%%%%%%%%%%%%%%%%%

\subsection{Nonsplitting MASAs of ${\mathbf e(p_0+k_0,q_0+k_0)}$ related to non-free-rowed MANSs}

The general study of non-free rowed MASAs of $o\pq$ is less well developed. Many different 
series of MASAs of $o\pq$ exist.
We will consider only two of them, which we denote $A(2k+1,0)$ and $A(2k+1,1)$, 
by analogy with series of non-free-rowed MANSs of $o(n,\mC)$ \cite{verc}.

\noindent
1. The series $A(2k+1,0)$ of $o\pq$ is represented by the matrix set
\bea
X =\left( \begin{array}{ccccccc} 0 & a_1 & 0 & a_2 & \ldots & a_k & 0 \\
                           & \ddots & \ddots & \ddots & \ddots & & a_k \\
                             &  & \ddots & \ddots & \ddots & \ddots & \vdots \\
                              & & &  \ddots & \ddots & \ddots & a_2 \\
                                & & & &  \ddots & \ddots & 0 \\
                               & & & & &  \ddots & a_1 \\
                               & & & & & & 0
\end{array} \right),  \label{splserie1}\\
K=F_{2k+1}=\left( \begin{array}{ccccc} & & & & \epsilon \\
                               & & & -\epsilon & \\
                                & & \invddots & &  \\
                                 & -\epsilon & & & \\
                                \epsilon & & & &
\end{array} \right),  
\eea
where all ${a_i's}$ are free. \\
Thus for $\epsilon=1$ we have $M \subset \left\{ \begin{array}{cc} o(k+1,k) & 
{\rm for} \, \, k \, \,{\rm even}   \\ 
 o(k,k+1) & {\rm for} \, \, k \, \, {\rm odd} 
                              \end{array} \right. $ \\
and for $\epsilon=-1$ we have $M \subset  \left\{ \begin{array}{cc}  o(k+1,k) & 
{\rm for} \, \, k \, \, {\rm odd}\\ 
o(k,k+1) & {\rm for} \, \, k \, \, {\rm even.} 
                              \end{array} \right.  $   \\

The splitting MASA of $e\pq$ for this series (in accordance with Theorem 4.1) is written as follows:
\bea
X_e =\left( \begin{array}{cccccccc} 0 & a_1 & 0 & a_2 & \ldots & a_k & 0 & \alpha \\
                           & \ddots & \ddots & \ddots & \ddots & & a_k & 0\\
                             &  & \ddots & \ddots & \ddots & \ddots & \vdots & \vdots \\
                              & & &  \ddots & \ddots & \ddots & a_2  & 0 \\
                                & & & &  \ddots & \ddots & 0 & 0 \\
                               & & & & &  \ddots & a_1 & 0 \\
                               & & & & & & 0 & 0 \\
                               & & & & & & & 0 
\end{array} \right). \label{seria1}
\eea

\bth
Every nonsplitting MASA of $e\pq$ corresponding to the splitting MASA (\ref{seria1}) is $E\pq$ conjugate to the 
following one
\bea
X_e =\left( \begin{array}{ccccccccc} 0 & a_1 & 0 & a_2 & \ldots & \ldots & a_k & 0 & \alpha \\
                           & \ddots & \ddots & \ddots & \ddots & & & a_k & 0\\
                             &  & \ddots & \ddots & \ddots & \ddots & & 0 & a_k \\
                                     & & & \ddots & \ddots & \ddots & \ddots  & \vdots & \vdots \\
                              & & &  & \ddots & \ddots & \ddots  & a_2  & 0 \\
                                & & & & &  \ddots & \ddots & 0 & a_2\\
                               & & & & & & \ddots & a_1 & 0\\
                               & & & & & & & 0 & a_1 \\
                               & & & & &  & & & 0 
\end{array} \right), \, \,
K_e=\left( \begin{array}{cc} F_{2k+1} & \\
                                    & 0
            \end{array} \right), \label{nonseria1}
\eea
where all entries in $X_e$ are free.

\enth

{\it Proof:} We will construct a nonsplitting MASA from the splitting one (\ref{seria1})

\bea
X'_e =\left( \begin{array}{ccccccccc} 0 & a_1 & 0 & a_2 & \ldots & & a_k & 0 & \alpha \\
                           & 0 & a_1 & 0 & a_2 & \ldots & & a_k & \beta_2\\
                             &  & \ddots & \ddots & \ddots & \ddots  & & 0 & \beta_3 \\
                               & & & \ddots & \ddots & \ddots & \ddots & \vdots & \vdots \\
                              & & & & \ddots & \ddots & \ddots & a_2  & \beta_{2k-2} \\
                               & & & & &  \ddots & \ddots & 0 & \beta_{2k-1} \\
                               & & & & & &  \ddots & a_1 & \beta_{2k}\\
                               & & & & & & & 0 & \beta_{2k+1} \\
                               & & & & & & & & 0 
\end{array} \right), \label{ll}
\eea
where $\beta 's$ are linearly dependent on $a_i 's$. Before imposing commutation 
relations we will remove the coboundaries. 

Consider one element of the algebra (\ref{ll})
\bea
A_1= \left( \begin{array}{ccccccccc}  0 & 1 & 0 & 0 & \ldots & & \ldots  & 0 & 0 \\
                           & 0 & 1 & 0 & \ldots & & \ldots & 0 & \alpha_{1,2} \\
                             &  & \ddots & \ddots & \ddots &  & &  0 & \alpha_{1,3} \\
                               & & & \ddots & \ddots & \ddots &  & \vdots & \vdots \\
                              & & & &  \ddots & \ddots & \ddots & 0  & \alpha_{1,2k-2} \\
                                & & &&  &  \ddots & \ddots & 0 & \alpha_{1,2k-1} \\
                               & & & & & &  \ddots & 1 & \alpha_{1,2k}\\
                               & & & & & & & 0 & \alpha_{1,2k+1} \\
                               & & & & & & & & 0 
\end{array} \right), \label{a1}
\eea
where $\alpha_{1,l}, l=2, \ldots, 2k+1$ represent the translations. We note that 
$\alpha_{1,l} \ldots \alpha_{1,2k}$ correspond to coboundaries and can be eliminated 
by conjugation by the translation group. Thus only $\alpha_{1,2k+1}$ is left in $A_1$.

Now consider an element $A_i$ of algebra (\ref{ll}), obtained by setting  
$a_i=\delta_{ij}$, $j \geq2$
\bea
A_i= \left( \begin{array}{ccccccccc}  0 & 0 & 0 & 1 & \ldots & & \ldots & 0 & 0 \\
                           & 0 & 0 & 0 & 1 & & \ldots & 0 & \alpha_{i,2} \\
                             &  & \ddots & \ddots & \ddots & \ddots &  & 0 & \alpha_{i,3} \\
                              & & & \ddots & \ddots & \ddots & \ddots & \vdots & \vdots \\
                              & & & &  \ddots & \ddots & \ddots & 1  & \alpha_{i,2k-2} \\
                                & & & &  &  \ddots & \ddots & 0 & \alpha_{i,2k-1} \\
                               & & & & & &  \ddots & 0 & \alpha_{i,2k}\\
                               & & & & & & & 0 & \alpha_{i,2k+1} \\
                               & & & & & & & & 0 
\end{array} \right) \label{ai}.
\eea
Commuting $A_1$ with all $A_i, i=2, \ldots k$  we obtain that  $\alpha_{j,2k-2j+3} = \alpha_{1,2k+1}$, \linebreak
$j=2, \ldots,  k$ and all other $\alpha_{i,j}$ have to be zero.

Using the normalizer $G$ of the form
\bea
G=\left( g_k^k, \ldots , g_k^2, g_k, 1, g_k^{-1}, \ldots, g_k^{-k} \right)
\eea
we can normalize  $\alpha_{1,2k+1}$ to $\alpha_{1,2k+1}=1$. This leads to the 
MASA (\ref{nonseria1}) and completes the proof of Theorem 5.2.
$\hfill\Box$

\noindent
2. The series $A(2k+1,1)$ of $o\pq$ is represented by the following matrix set:
\bea
X =\left( \begin{array}{cccccccc} 0 & a_1 & 0 & a_2 & \ldots & a_k & 0 & b\\
                           & \ddots & \ddots & \ddots & \ddots & & a_k & 0\\
                             &  & \ddots & \ddots & \ddots & \ddots & \vdots & \vdots \\
                              & & &  \ddots & \ddots & \ddots & a_2 & 0\\
                                & & & &  \ddots & \ddots & 0 & 0\\
                               & & & & &  \ddots & a_1 & 0 \\
                               & & & & & & 0 & 0\\
                               & & & & & 0 & -\epsilon b & 0 
\end{array} \right), \qquad
K=\left( \begin{array}{cc} F_{2k+1} & \\ & 1 \end{array} \right),
\label{serie2}
\eea
where all ${a_i's}$ and $b$ are free. The corresponding metric is
\bea
K=\left( \begin{array}{cc} F_{2k+1} & \\ & 1 \end{array} \right) =
\left( \begin{array}{cccccc} & & & & \epsilon & 0\\
                               & & & -\epsilon & & 0\\
                                & & \invddots & & & \vdots \\
                                 & -\epsilon & & & & 0\\
                                \epsilon & & & & & 0 \\
                                0 & 0 & \ldots & \ldots & 0 & 1 
\end{array} \right) \label{Kseria2}
\eea
Thus for $\epsilon=1$ we have $M \subset \left\{ \begin{array}{ll} o(k+1,k+1) & {\rm for} \, \, k \, \,{\rm odd}\\ 
                                                         o(k+2,k) & {\rm for} \, \, k \, \, {\rm even} 
                              \end{array} \right. $  \\
and for $\epsilon=-1$ we have $M \subset \left\{ \begin{array}{ll} o(k+1,k+1) & {\rm for} \, \, k \, \, {\rm even}\\ 
                                                         o(k+2,k) & {\rm for} \, \, k \, \, {\rm odd.} 
                              \end{array} \right. $  \\

\bth
Every nonsplitting MASA corresponding to the splitting MASA (\ref{serie2}) is $E\pq$ conjugated to the
MASA of the form
\bea
X_e =\left( \begin{array}{cccccccccc} 0 & a_1 & 0 & a_2 & \ldots & & a_k & 0 & b & \alpha\\
                           & \ddots & \ddots & \ddots & \ddots & & & a_k & 0 & \lambda b\\
                             &  & \ddots & \ddots & \ddots & \ddots & & 0 & 0 & 0\\
                                & & & \ddots & \ddots & \ddots & \ddots & \vdots & \vdots & \vdots \\
                              & & & & \ddots & \ddots & \ddots & a_2 & 0 & 0 \\
                                & & & &  &  \ddots & \ddots & 0 & 0 & 0\\
                               & & & & & &  \ddots & a_1 & 0 & 0 \\
                               & & & & & & & 0 & 0 & 0 \\
                               & & & & & & & -\epsilon b & 0 & \lambda a_1 + \mu b \\
                               & & & & & & & 0 & 0 & 0_1
\end{array} \right) \label{nonsplserie2}
\eea
with the metric as in (\ref{Kseria2}). The entries $a_i, b$ and $\alpha$ are free. Parameters
$\lambda$ and $\mu$ are one of the following sets:
\bea
(\lambda, \mu)= \left\{ \begin{array}{l} (0,1)\\ (0,-1) \\ (1, \mu), \, \, \mu \in \mR .
\end{array} \right. 
\eea
\enth
{\it Proof:} The proof is similar to that of Theorem 5.2 and we omit it here.
$\hfill\Box$

%%%%%%%%%%%%%%%%%%%%%%%%%%%%%%%%%%%%%%%%%%%%%%%%%%%%%%%%%%%%%%%%%%%
%%%%%%%%%%%%%%%%%%%%%%%%%%%%%%%%%%%%%%%%%%%%%%%%%%%%%%%%%%%%%%%%%%%
%%%%%%%%%%%%%%%%%%%%%                         %%%%%%%%%%%%%%%%%%%%%
%%%%%%%%%%%%%%%%%%%%%%%%%%%%%%%%%%%%%%%%%%%%%%%%%%%%%%%%%%%%%%%%%%%
%%%%%%%%%%%%%%%%%%%%%%%%%%%%%%%%%%%%%%%%%%%%%%%%%%%%%%%%%%%%%%%%%%%
%%%%%%%%%%%%%%%%%%%%%%%%%%%%%%%%%%%%%%%%%%%%%%%%%%%%%%%%%%%%%%%%%%%
%%%%%%%%%%%%%%%%%%%%%%%%%%%%%%%%%%%%%%%%%%%%%%%%%%%%%%%%%%%%%%%%%%%
%%%%%%%%%%%%%%%%%%%%%                         %%%%%%%%%%%%%%%%%%%%%
%%%%%%%%%%%%%%%%%%%%%%%%%%%%%%%%%%%%%%%%%%%%%%%%%%%%%%%%%%%%%%%%%%%
%%%%%%%%%%%%%%%%%%%%%%%%%%%%%%%%%%%%%%%%%%%%%%%%%%%%%%%%%%%%%%%%%%%
\section{Decomposition properties of MASAs of ${\mathbf e\pq}$}
\setcounter{equation}{0}

The results of Sections 4 and 5 can be formulated in terms of a decomposition of
the underlying pseudoeuclidean space $S\pq$. Both splitting and nonsplitting MASAs 
have been represented by matrix sets $\{X_e,K_e\}$ as in eq.(\ref{delta0}), (\ref{Ke1}).
We shall call a MASA of $e\pq$ {\it decomposable} if the metric $K_e$ in  (\ref{Ke1})
consists of 2 or more blocks. The projection of such a MASA onto the $o\pq$ subalgebra is 
then an orthogonally decomposable MASA of $o\pq$. Let $M_e\pq$ be a decomposable MASA 
of $e\pq$.
The space $S\pq$ then splits into a direct sum of subspaces
\bea
S\pq =\bigoplus_{i=1}^{l} S(p_i,q_i), \qquad \sum_{i=1}^l p_i =p, \qquad \sum_{i=1}^l q_i =q
\eea
and each indecomposable component of the decomposable MASA of $e\pq$ acts independently 
in one of the spaces $S(p_i,q_i)$. We shall write
\bea
M_e\pq =\bigoplus_{i=1}^{l} M_e(p_i,q_i) \label{pattern}.
\eea
Each individual indecomposable MASA $M_e(p_i,q_i) \subset e(p_i,q_i)$ can then be considered 
separately.

Consider the matrix set $\{X_e,K_e\}$, $X_e$ given by eq.(\ref{delta0}), $K_e$ as in eq.(\ref{Ke1}),
where each block is indecomposable. The blocks to be considered consist of a block on the diagonal in $X_e$, 
plus an entry from the right hand column in $X_e$.

The following types of indecomposable MASAs $M_e(p_i,q_i) \subset e(p_i,q_i)$ exist.
\begin{itemize}
\item
$dimS=1$. The MASAs are pure positive or negative length translations.
\bea
&M_e(1,0) = & \left\{  \left( \begin{array}{cc} 0 & x \\ 0 & 0 \end{array} \right), \quad  x \in \mR, \quad K_e=\left( \begin{array}{cc} 1 & 0 \\ 0 & 0 \end{array} \right) \right\} \\
&M_e(0,1)  = & \left\{ \left( \begin{array}{cc} 0 & y \\ 0 & 0 \end{array} \right), \quad y \in \mR, \quad K_e=\left( \begin{array}{cc} -1 & 0 \\ 0 & 0 \end{array} \right) \right\}.
\eea
A MASA $M_e(p,q)$ of $e\pq$ contains $\kp$ of the first ones and $\km$ of the second.
\item
$dimS=2$. The MASAs are $o(2)$ rotations in a $(++)$, or $(--)$ type subspace, or $o(1,1)$ pseudorotations in a $(+-)$ space:
\bea
& M_e(2,0) = & \left\{ \left( \begin{array}{ccc} 0 & x & 0\\ -x & 0 & 0 \\ 0 & 0 & 0 \end{array} \right) , \quad
K_e= \left( \begin{array}{cc} I_2 & 0 \\ 0 & 0 \end{array} \right) \right\} \\
& M_e(0,2) = & \left\{ \left( \begin{array}{ccc} 0 & x & 0\\ -x & 0 & 0 \\ 0 & 0 & 0 \end{array} \right) , \quad
K_e= \left( \begin{array}{cc} -I_2 & 0 \\ 0 & 0 \end{array} \right) \right\} \\
& M_e(1,1) = & \left\{ \left( \begin{array}{ccc} a & 0 & 0\\ 0 & -a & 0 \\ 0 & 0 & 0 \end{array} \right) , \quad
K_e= \left( \begin{array}{ccc} 0&1 & 0 \\ 1&0 & 0 \\0 & 0 & 0 \end{array} \right) \right\}. 
\eea
\item
$dimS=k \geq 3$. There are two possible types of indecomposable MASAs of $e\pq$ for $p+q \geq 3$. Both of them have $\kp=\km=0$ (no nonisotropic translations).
\begin{zoznamrom}
\item $M_e\pq$ contains $\kz$ isotropic translations with $\kz \geq 1$. The projection of $M_e\pq$ onto $o\pq$ is then a MANS of $o\pq$ with Kravchuk signature $(\kz, p+q-2\kz,\kz)$. The MANS can be free-rowed or non-free-rowed. The MASA of $e\pq$ can be splitting , or nonsplitting. Such MASAs exist for any $p+q \geq 3$, $min(p,q) \geq 1$. They were treated in Sections 4 and 5.
\item $M\pq$ is an orthogonally indecomposable MASA of $o\pq$ that is not a MANS. It gives  rise to a splitting MASA of $e\pq$ which contains no translations ($\kz=0$). As reviewed in Section 3 such MASAs of  $o\pq$ exist only for $p+q$ even.

\end{zoznamrom}
\end{itemize}

%%%%%%%%%%%%%%%%%%%%%%%%%%%%%%%%%%%%%%%%%%%%%%%%%%%%%%%%%%%%%%%%%%%
%%%%%%%%%%%%%%%%%%%%%%%%%%%%%%%%%%%%%%%%%%%%%%%%%%%%%%%%%%%%%%%%%%%
%%%%%%%%%%%%%%%%%%%%%                         %%%%%%%%%%%%%%%%%%%%%
%%%%%%%%%%%%%%%%%%%%%%%%%%%%%%%%%%%%%%%%%%%%%%%%%%%%%%%%%%%%%%%%%%%
%%%%%%%%%%%%%%%%%%%%%%%%%%%%%%%%%%%%%%%%%%%%%%%%%%%%%%%%%%%%%%%%%%%
%%%%%%%%%%%%%%%%%%%%%%%%%%%%%%%%%%%%%%%%%%%%%%%%%%%%%%%%%%%%%%%%%%%
%%%%%%%%%%%%%%%%%%%%%%%%%%%%%%%%%%%%%%%%%%%%%%%%%%%%%%%%%%%%%%%%%%%
%%%%%%%%%%%%%%%%%%%%%                         %%%%%%%%%%%%%%%%%%%%%
%%%%%%%%%%%%%%%%%%%%%%%%%%%%%%%%%%%%%%%%%%%%%%%%%%%%%%%%%%%%%%%%%%%
%%%%%%%%%%%%%%%%%%%%%%%%%%%%%%%%%%%%%%%%%%%%%%%%%%%%%%%%%%%%%%%%%%%

\section{A special case: MASAs of e(p,2)}
\setcounter{abc}{0}

The case $q=2$, like $q=1$ and $q=0$, presented earlier \cite{conf} is simpler than that 
of $q \geq 3$. All MASAs can be presented explicitly, in particular those involving non-free-rowed MANS of $o(p,2)$.

The possible decomposition patterns (\ref{pattern}) for MASAs of $e(p,2)$ are
\bea
M_e(p,2) & = & M_e(p_1,2) \oplus l_+M_e(2,0) +\kp M_e(1,0) \\
&  & p_1=1, \quad or \quad p_1 \geq 2, \qquad p_1+2 l_+ + \kp =p  \nonumber \\
M_e(p,2) & = & M_e(p_1,1) \oplus M_e(p_2,1)\oplus l_+M_e(2,0) +\kp M_e(1,0) \\
& & p_1+p_2+2 l_+ + \kp =p  \nonumber \\
M_e(p,2) & = & M_e(0,2) \oplus l_+M_e(2,0) +\kp M_e(1,0) \\
&  & 2 l_+ + \kp =p.  \nonumber
\eea
The algebras $M_e(2,0)$, $M_e(0,2)$ and $M_e(1,0)$ are already abelian (and one dimensional) as
are $M_e(0,1)$ and $M_e(1,1)$. The MASAs $M_e(p,1)$ of $e(p,1)$, $p \geq 2$ were studied in our earlier article \cite{conf}.

Thus, we  need  to treat only indecomposable MASAs of $e(p,2)$. As was stated in Section 6 for general $e\pq$, two cases 
arise, namely $\kz=0$ and $1 \leq \kz \leq min(p,q)$, where $\kz$ is the number of linearly independent translation generators present. 

\noindent
1. $\kz=0$ \\
Then $M(p,2)$ is an orthogonally indecomposable MASA of $o(p,2)$ that is not a MANS. These exist only when $p$ is even ($p \geq 2$). 

For $p=2$ three inequivalent OID MASAs that are not MANS exist and the corresponding splitting MASAs of $e(p,2)$ are given by the following matrix sets:
\begin{zoznamrom}
\item $M(2,2)$ is AOID but D
\bea X_{e}=\left( \begin{array}{ccccc} 
                               a & b & \ & \ & 0\\
                               0 & a & \ & \  & 0\\
                               \ & \ & -a & 0  & 0\\
                                   \ & \ & -b & -a & 0 \\
                                   & & &  & 0_1
                       \end{array} \right), \quad
K_{e}=\left( \begin{array}{ccc}   \ & I_{2} & \  \\
                                         I_{2} & \ &  \  \\
                                 \ & \ &  0_{1}
                      \end{array} \right) \protect\label{p2spl:k00Ci}
\eea
\item $M(2,2)$ is AOID, ID but NAID  
\bea X_{e}=\left( \begin{array}{ccccc} 
                               0 & a & 0 & b & 0 \\
                              -a & 0 & -b & 0 & 0 \\
                             \ & \ & 0 & a & 0 \\
                              \ & \ & -a & 0  & 0 \\
                                 & & & & 0_1
                       \end{array} \right)
\protect\label{p2spl:k00Cii}
\eea
with $K_{e}$ same as in {\em i)}. 
\item $M(2,2)$ is NAOID but D 
\bea X_{e}=\left( \begin{array}{ccccc}
                               a & b & \ & \ & 0 \\
                              -b & a & \ & \  & 0\\
                               \ & \ & -a & b & 0 \\
                               \ & \ & -b & -a & 0 \\
                                  & & & & 0_1
                       \end{array} \right)
\protect\label{p2spl:k00Ciii}
\eea
with $K_{e}$ same as in {\em i)}.
\end{zoznamrom}

For $p=2l$, $l\geq2$ we have just one OID MASA of $o(p,2)$ (NAOID, ID but NAID), namely
$M=RQ \oplus {\rm MANS \, \, of \, \,} su(l,1)$. The corresponding splitting MASA of $e(p,2)$ is
represented as following matrix set
\begin{eqnarray*} X_{e}=\left( \begin{array}{cccccccccc}
       0 & b & a_{1} & 0 & \ldots & a_{l-1} & 0 & 0 & c & 0 \\
      -b & 0 & 0 & a_{1} & \ldots & 0 & a_{l-1} & -c & 0 & 0  \\
       \ & \ & 0 & b & \ & \ & \ & -a_{1} & 0  & 0\\ 
       \ & \ & -b & 0 &  \ & \ & \ & 0 & -a_{1} & 0  \\
       \ & \ & \ & \ & \ddots & \ & \ & \vdots & \vdots & \vdots \\
       \ & \ & \ & \ & \ & 0 & b & -a_{l-1} & 0 & \vdots \\
       \ & \ & \ & \ & \ & -b & 0 & 0 &  -a_{l-1} & 0 \\
       \ & \ & \ & \ & \ & \ & \ & 0 & b & 0 \\
       \ & \ & \ & \ & \ & \ & \ & -b & 0 & 0 \\
          &   &   &   &  &   &   &    &    & 0_1
                               \end{array} \right),
 \protect\label{p2spl:k00D}
\end{eqnarray*} 
\bea K_{e}=\left( \begin{array}{cccc} \ & \ & I_{2} & \  \\
                                       \ & I_{2l-2} & \ & \  \\
                                        I_{2} & \ & \ & \  \\
                                        \ & \ & \ & 0_{1} 
                   \end{array} \right). 
\eea

\noindent
2. $\kz=1$ \\
The projection of $M_e(p,2)$ onto $o(p,2)$ will be a MANS of $o(p,2)$ with Kravchuk 
signature $( 1 \, \, p \, \, 1)$. This MANS can be free-rowed, or non-free-rowed, so we 
obtain two splitting MASAs of $e(p,2)$ represented, respectively, by
\begin{zoznamrom}
\item  free-rowed
\bea X_{e}=\left( \begin{array}{cccc} 
                             0 & \alpha & 0 & z \\
                             0 & 0 & -K_{0} \alpha^{T} & 0  \\
                             0 & 0 & 0  & 0 \\
                               &   &  &  0_1
                          \end{array} \right), \quad
K_{e}=\left( \begin{array}{ccc} 0 & 0 & 1    \\
                                 0 & K_0 & 0  \\
                                1 & 0 & 0   
             \end{array} \right), \protect\label{p2spl:k01Bfr}
\eea
where $ K_{0}$ has  signature $(p-1,1)$, $\alpha \in {\mR}^{1 \times p}$, 
$1 \leq p$ 
\item non-free rowed  
\bea
\begin{array}{l} X_{e}=\left( \begin{array}{ccccccc} 
                              0 & a & \alpha & 0 & b & 0 & z \\
                              \ & 0 & 0 & a & 0 & -b  & 0\\
                              \ & \ & 0 & 0 & 0 & -\alpha^{T} & \vdots\\
                               \ & \ & \ & 0 & -a & 0 & 0 \\
                              \ & \ & \ & \ & 0 & -a & 0 \\
                               & \ & \ & \ & \ & 0 & 0 \\
                                   & & &  &  & & 0_1
                   \end{array} \right), \protect\label{p2spl:k01Bnfr} \\
K_{e}=\left( \begin{array}{cccccc} 
                                \ & \ & \ & \ & 1 & 0 \\
                                \ & \ & \ & 1 & 0  & 0\\
                                 \ & \ & I_{\nu +1} & 0 & 0  & \vdots \\
                                 \ & 1 & 0 & 0 & 0 & 0  \\
                                 1 & 0 & 0 & 0 & 0 & 0 \\
                                   & & & & & 0_1
                      \end{array} \right)  
\end{array}
\eea 
$\alpha \in {\mR}^{1 \times \nu}$, $1 \leq \nu $ and $\nu = p-3$. 
\end{zoznamrom}

The MASA (\ref{p2spl:k01Bfr}) gives rise to three different nonsplitting MASAs for $p \geq2$
which can be expressed as 
\bea X_{e}=\left( \begin{array}{cccc} 0 & \alpha & 0 & z \\
                             0 & 0 & -K_{0}\alpha^{T} & BK_{0}\alpha^{T} \\
                                     0 & 0 & 0 & 0 \\
                                       0 & 0 & 0 & 0 
                  \end{array} \right), \quad
 K_{e}=\left( \begin{array}{cccc} \ & \ & 1 & \ \\
                                      \ & K_{0} & \ & \ \\
                                     1 & \ & \ & \ \\
                                      \ & \ & \ & 0_{1}
                    \end{array} \right).
\eea
$K_{0}$ is the same as in (\ref{p2spl:k01Bfr}) and B satisfies the condition $ BK_{0}=K_{0}B^{T}$, {\it i.e.} B is an element of the Jordan algebra $jo(p-1,1)$.
A classification of the elements Jordan algebras was performed in the paper by Djokovic et al \cite{dj} and the couple $\{B,K_{0}\}$ can have one of  the three different following forms (keeping in mind the signature of $K_0$):
\setcounter{abc}{0}
\begin{zoznamrom}
\item \bea B=\left( \begin{array}{cc} a & \ \\
                                         \ & B_{0} 
                       \end{array} \right), \quad
K_{0}=\left( \begin{array}{cc} -1 & \ \\
                                     \ & I  
                   \end{array} \right) 
\eea
\newline
\item \bea B=\left( \begin{array}{ccc} a & 0 & \ \\
                                         1 & a & \ \\
                                         \ & \ & B_{0} 
                        \end{array} \right), \quad 
K_{0}=\left( \begin{array}{ccc} 0 & 1 & \ \\
                                     1 & 0 & \ \\
                                     \ & \ & I 
                   \end{array} \right) 
\eea
\newline
\item \bea B=\left( \begin{array}{cccc} a & 0 & 0 & \ \\ 
                                            1 & a & 0 & \ \\
                                            0 & 1 & a & \ \\
                                            \ & \ & \ &  B_{0}
                       \end{array} \right), \quad 
K_{0}=\left( \begin{array}{cccc} 0 & 0 & 1 & \ \\ 
                                      0 & 1 & 0 & \ \\  
                                      1 & 0 & 0 & \ \\
                                      \ & \ & \ & I 
                  \end{array} \right), 
\eea
where $  B_{0}$ is a diagonal matrix.
\end{zoznamrom}

For $p=1$ the nonsplitting MASA corresponding to eq.(\ref{p2spl:k01Bfr}) is 
\bea X_{e}=\left( \begin{array}{cccc} 0 & a & 0 & z \\
                                      0 & 0 & -a & 0 \\
                                      0 & 0 & 0 & a \\
                                      0 & 0 & 0 & 0
                   \end{array} \right), \quad
K_{e}=\left( \begin{array}{cccc} 0 & 0 & 1 & 0 \\
                                      0 & 1 & 0 & 0 \\
                                     1 & 0 & 0 & 0 \\
                                     0 & 0 & 0 & 0_{1}
                \end{array} \right) 
\eea

\noindent
The MASA (\ref{p2spl:k01Bnfr}) for $\nu \geq 2$ gives rise to one type of nonsplitting MASA that 
can be represented as
\bea X_{e}=\left( \begin{array}{ccccccc}
             0 & a & \alpha & 0 & b & 0 &  z \\ 
          \ & 0 & 0 & a & 0 & -b &   \alpha\rho^{T} \\
     \ & \ & 0 & 0 & 0 & -\alpha^{T} &  a\rho^T + \Lambda \alpha^{T} \\
              \ & \ & \ & 0 & -a & 0 &  0 \\
               \ & \ & \ & \ & 0 & -a &  0 \\
               0 & \ & \ & \ & \ & 0 &  0 \\
              \ & \ & \ & \ & \ & \ & 0_{1} 
\end{array} \right) 
\eea
with $\Lambda=\Lambda^T$.
Using the normalizer $G=diag(g,g_1,G_2,g_3,1/g_1,g,1)$, $G_2 \in \mR^{\nu \times \nu}$, $g,g_1,g_3 \in \mR$,  satisfying $G_2G_2^T=I_{\nu}$, $g^2=g_3^2=1$ we can transform $\Lambda,\rho)$ into
\beq
\Lambda'={1 \over g} G_2 \Lambda G_2^T, \qquad \rho'={1 \over g_1g_3} G_2 \rho.
\eeq
We can use $G_2$ either to diagonalize $\Lambda$, or to rotate $\rho$ into {\it e.g.} $\rho=(\rho_1, 0, \ldots, 0)$.

\noindent
3. $\kz =2$ \\
The projection of $M_e(p,2)$ onto $o(p,2)$ is a free-rowed MANS with Kravchuk signature $(2 \, \, p\!-\!2 \, \, 2)$. The corresponding splitting MASA of $e(p,2)$ is given in Theorem 5.1 with $q=\kz=2$ and $K_{p_0,q_0}=I_{p-2}$. In this case $Q_2$ can be chosen as $Q_2=diag(1,q_2, \ldots , q_{\mu})$, $q_{1}=1 \geq |q_{2}| \geq \ldots \geq |q_{\mu}|$.
\noindent
This MASA in turn gives rise to the following non-splitting MASAs.
\bea X_{e}=\left( \begin{array}{cccccc} 0 & 0 & \alpha & 0 & y & z_{1} \\
                               0 & 0 & \alpha Q & -y & 0 & z_{2} \\
                                0 & 0 & 0 & -\alpha^{T} & -Qa^{T} & \Lambda^T \alpha^T \\
                                 0 & 0 & 0 & 0 & 0 & 0 \\
                                 0 & 0 & 0 & 0 & 0 & 0 \\
                                 0 & 0 & 0 & 0 & 0 & 0_1 
                  \end{array} \right). 
\eea
Here $\Lambda$ is a diagonal matrix, $Tr\Lambda =0$ and $K_e$ is same as 
in eq.(\ref{Zt}).

%%%%%%%%%%%%%%%%%%%%%%%%%%%%%%%%%%%%%%%%%%%%%%%%%%%%%%%%%%%%%%%%%%%
%%%%%%%%%%%%%%%%%%%%%%%%%%%%%%%%%%%%%%%%%%%%%%%%%%%%%%%%%%%%%%%%%%%
%%%%%%%%%%%%%%%%%%%%%                         %%%%%%%%%%%%%%%%%%%%%
%%%%%%%%%%%%%%%%%%%%%%%%%%%%%%%%%%%%%%%%%%%%%%%%%%%%%%%%%%%%%%%%%%%
%%%%%%%%%%%%%%%%%%%%%%%%%%%%%%%%%%%%%%%%%%%%%%%%%%%%%%%%%%%%%%%%%%%

%%%%%%%%%%%%%%%%%%%%%%%%%%%%%%%%%%%%%%%%%%%%%%%%%%%%%%%%%%%%%%%%%%%
%%%%%%%%%%%%%%%%%%%%%%%%%%%%%%%%%%%%%%%%%%%%%%%%%%%%%%%%%%%%%%%%%%%
%%%%%%%%%%%%%%%%%%%%%                         %%%%%%%%%%%%%%%%%%%%%
%%%%%%%%%%%%%%%%%%%%%%%%%%%%%%%%%%%%%%%%%%%%%%%%%%%%%%%%%%%%%%%%%%%
%%%%%%%%%%%%%%%%%%%%%%%%%%%%%%%%%%%%%%%%%%%%%%%%%%%%%%%%%%%%%%%%%%%

\section {Conclusions}
The main conclusion is that we have presented guidelines for constructing all 
MASAs of $e\pq$ for any fixed values of $p$ and $q$. Some of the results are 
entirely explicit, such as Theorem 4.1 describing all splitting MASAs of $e\pq$, 
and Theorem 5.1 presenting nonsplitting MASAs containing a free-rowed MANS of 
$o(p_0+\kz,q_0+\kz) \subset o(p,q)$. The results on MASAs of $e\pq$ involving 
non-free-rowed MANS of $o(p_0+\kz,q_0+\kz)$ are less complete and amount to 
specific examples (see Theorems 5.2 and 5.3). The decomposition results of 
Section 6 allow us to restrict all considerations to indecomposable MASAs of 
$e\pq$, both splitting and non-splitting ones. The results for $e(p,2)$ presented 
in Section 7 are complete and explicit, like those given earlier for $e(p,0)$ 
and $e(p,1)$ \cite{conf}. In particular we have constructed all MASAs related 
to non-free-rowed MANSs.

Work concerning the application of MASAs of $e\pq$ is in progress. In particular,
we use  MASAs of $e\pq$ to construct  the coordinate 
systems in which certain partial differential equations (Laplace-Beltrami, 
Hamilton-Jacobi) 
allow the separation of variables.

\bigskip
\bigskip

\noindent
{\large \it Acknowledgement}

The research of P.W. was partially supported by research grants from NSERC of 
Canada and FCAR du Qu\'ebec.

%%%%%%%%%%%%%%%%%%%%%%%%%%%%%%%%%%%%%%%%%%%%%%%%%%%%%%%%%%%%%%%%%%%
%%%%%%%%%%%%%%%%%%%%%%%%%%%%%%%%%%%%%%%%%%%%%%%%%%%%%%%%%%%%%%%%%%%
%%%%%%%%%%%%%%%%%%%%%                         %%%%%%%%%%%%%%%%%%%%%
%%%%%%%%%%%%%%%%%%%%%%%%%%%%%%%%%%%%%%%%%%%%%%%%%%%%%%%%%%%%%%%%%%%
%%%%%%%%%%%%%%%%%%%%%%%%%%%%%%%%%%%%%%%%%%%%%%%%%%%%%%%%%%%%%%%%%%%

\end {document}